\documentclass[12pt]{article}
\usepackage{graphicx}
\usepackage{mathptmx}
\usepackage{latexsym,amsmath,amssymb,amsfonts,amsthm}

\newcommand{\ti}[1]{\mbox{\tiny $#1$}}

\newcommand{\ft}[1]{\mbox{\footnotesize $#1$}}
\newcommand{\sm}[1]{\mbox{\small $#1$}}
\newcommand{\la}[1]{\mbox{\large $#1$}}
\newcommand{\La}[1]{\mbox{\Large $#1$}}
\newcommand{\LA}[1]{\mbox{\LARGE $#1$}}

\newfont{\lie}{eufm10 at 12pt}
\newfont{\field}{msbm10 at 11pt}

\newtheorem{theorem}{Theorem}[section]
\newtheorem{lemma}{Lemma}[section]
\newtheorem{corollary}{Corollary}[section]
\newtheorem{proposition}{Proposition}[section]
\newtheorem{remark}{Remark}[section]

\newtheorem{definition}{Definition}[section]

\newtheorem{example}{Example}[section]
\begin{document}
\title{\normalsize  \bf   STABILITY OF
 SUBMANIFOLDS WITH PARALLEL MEAN CURVATURE IN
CALIBRATED MANIFOLDS}
\markright{\sl 
\hfill STABILITY OF SUBMANIFOLDS WITH PARALLEL MEAN CURVATURE \hfill}
\author{ Isabel M.C.\ Salavessa}
\date{}
\protect\footnotetext{\!\!\!\!\!\!\!\!\!\!\!\!\! {\bf MSC 2000:}
Primary:  53C42; 53C38. Secondary: 58E12; 35J19; 47A75
Partially supported by
FCT through program
 PTDC/MAT/101007/2008.}
\maketitle ~~~\\[-10mm]
{\footnotesize Centro de F\'{\i}sica das Interac\c{c}\~{o}es
Fundamentais, Instituto Superior T\'{e}cnico, Technical University
of Lisbon, Edif\'{\i}cio Ci\^{e}ncia, Piso 3, Av.\ Rovisco Pais,
1049-001 Lisboa, Portugal;~ isabel.salavessa@ist.utl.pt}\\[5mm]
{\small {\bf Abstract:} On a 
Riemannian manifold $\bar{M}^{m+n}$ with an
$(m+1)$-calibration $\Omega$, we prove
 that an $m$-submanifold $M$ with
constant mean curvature $H$ and  calibrated extended
tangent space $\mathbb{R}H\oplus TM$ is a  critical point of the area 
functional for variations that preserve the enclosed $\Omega$-volume.
This recovers the case described by Barbosa, do Carmo and Eschenburg,
when $n=1$ and $\Omega$ is the volume element of $\bar{M}$.
To the second variation we  associate an
 $\Omega$-Jacobi operator and define $\Omega$-stability.
Under natural conditions, we show that the Euclidean $m$-spheres are the unique
$\Omega$-stable submanifolds of $\mathbb{R}^{m+n}$.
We study the $\Omega$-stability of  geodesic 
$m$-spheres of a fibred space form $M^{m+n}$ with
 totally geodesic $(m+1)$-dimensional fibres.
}
\section{Introduction} 
Immersed hypersurfaces  with constant mean curvature  of a Euclidean space 
 are known to be critical points of a variational problem, namely,
they are critical points of the $m$-area $A_D(t)$ for all
 variations $\phi_t:D\subset M^m\to \mathbb{R}^{m+1}$ of $\phi=\phi_0$
 fixing the boundary of
a compact domain $D$, and  
that leave a certain enclosed $(m+1)$-volume $V_D(t)$ invariant.
This volume can be given by $V_D(\phi)=\frac{1}{m+1}\int_D\langle
\phi,\nu\rangle dM$, where $\nu$ is the unit normal to $\phi$,
 and its modulo is the volume of the
cone over $\phi(D)$ with vertex at $0\in \mathbb{R}$ (see \cite{BdC}).
This property was generalized by Barbosa, do Carmo and
Eschenburg in \cite{BdCE} to hypersurfaces with constant
mean curvature $H$ immersed in a Riemannian manifold $\bar{M}^{m+1}$, 
by defining the volume of a variation $\bar{\phi}(t,p)=\phi_t(p)$ as
$V_D(t)=\int_{[0,t]\times D}\bar{\phi}^*d\bar{M}$.
A critical point of $A_D(t)$ for volume-preserving variations, i.e., 
 $V_D(t)=V_D(0)=0$, 
is just a  critical point of $J_D(t)=A_D(t)+mH_0V_D(t)$ for any variation
fixing the boundary,
where $H_0$ is the mean value of the mean curvature  $H$ of $\phi$, 
and it is characterized by having constant mean curvature $H_0$. 
Such a critical point is stable if
$A''_D(0)\geq 0$  for all volume preserving variations, or equivalently,
if $J''_D(0)\geq 0$ for all variations with vector variation $W$
satisfying $W^{\bot}=f\nu$, where 
  $f\in \mathcal{F}_D$ and  $\nu$ is the unit normal
of $M$. The class  $\mathcal{F}_D$ is given by the functions
$f:D\to \mathbb{R}$ such that $f=0$ on $\partial D$ and $\int_D fdM=0$.

Smooth solutions of  the isoperimetric problem which
seeks the least perimeter that encloses a given volume
 are stable hypersurfaces with constant mean curvature. 
So it is important to determine which hypersurfaces are stable.
If $\bar{M}$ is a space form, the  case
$M$ closed has been solved in \cite{BdCE}, concluding
 that $M$ must be a geodesic
sphere, and the   case $M$ complete, with $m=2$, has been partially solved 
by several authors (see e.g.\ the paper of Ritor\'{e} and Ros 
\cite{RiRo} and  references therein).
 
We ask if somehow we can extend these variational properties 
to higher codimension submanifolds of a Riemannian manifold
$\bar{M}$ of dimension $m+n$. We use a  pre-calibration
 $\Omega$ of rank $m+1$ on the ambient space $\bar{M}$ to define
the ``enclosed volume'' of a variation $\bar{\phi}: [0,\epsilon]
\times D\to \bar{M}$
as the $\Omega$-volume
$V_D(t)=\int_{[0,t]\times D}\bar{\phi}^*\Omega$. If $\bar{M}$ is
of dimension $m+1$ and $\Omega$ is its volume element, we recover
the case \cite{BdCE}. We  will assume $M$ has 
calibrated extended tangent space, that is, there exists a smooth
global unit normal $\nu$ such that $H=\|H\| \nu$ and
 $EM=\mathbb{R}\nu\oplus TM$ is a $\Omega$-calibrated vector bundle. 
This is a strong restriction, and  corresponds 
in some cases to
be able to  extend $M$  to
 a calibrated $(m+1)$-dimensional submanifold $M'$
such that   $TM'=EM$ along $M$. 
  Even in this  case,
 our approach differs from \cite{BdCE}, for we allow $\phi_t$ to take
values outside $M'$, and so
 our enclosed volume  at each time $t$
   may not correspond to the  enclosed volume
in $M'$ defined in \cite{BdCE}. 

Variational characterizations
of prescribed mean curvature was the subject
of earlier work of Gulliver \cite{Gu1,Gu2}, and Duzaar and Fuchs \cite{DF1,DF2}.
In \cite{Gu2} a stationary submanifold of a functional $A(D)+\int_D\alpha$, 
where $\alpha$ is an $m$-form such that $H=(d\alpha)^{\sharp}$, prescribes
the mean curvature $H$ of $M$ as an alternating $m$-tensor whose values are
orthogonal to each of its $m$ arguments. 
For a submanifold $M$ of a Euclidean space $\mathbb{R}^{m+n}$,
  Morgan in \cite{Mo}
defined a prescribed enclosed multi-volume, 
and proved  that $M$ is stationary for area for that 
prescribed multi-volume if and only if,
for some  $\xi\in \bigwedge_{m+1}\mathbb{R}^{m+n}$,
the mean curvature of $M$ satisfies  $H=\xi\lfloor \vec{S}$, where
$\vec{S}$ is the unit $m$-plane tangent to $M$.  This corresponds to
our condition of an $\Omega$-calibrated extended tangent space $EM$,
if  $\xi=\|H\|\Omega$, defining a  calibration
$\Omega$ (see Lemma 2.3).
Existence and regularity of such area-minimizing submanifolds 
(as rectifiable currents), with given boundary and multi-volume, 
are proved in \cite{Mo} under quite general
conditions.
  
We show that submanifolds with constant mean
 curvature are just the critical  points of $A_D(t)$ for variations that fix
the $\Omega$-volume, or equivalently, of $J_D(t)$
for any variation fixing the boundary $\partial D$. 
Furthermore, under certain conditions on
$\Omega$, it turns out these submanifolds 
have parallel mean curvature. 
We compute the second
variation of $J_D(t)$ and obtain 
$J_D''(0)=\int_D\bar{g}(\mathcal{J}'_{\ti{\Omega,D}}(W^{\bot}),W^{\bot})dM
=:I_{\Omega}(W^{\bot},W^{\bot})$, 
where $W=\frac{\partial \bar{\phi}}{\partial t}$  at $t=0$, and 
$$\mathcal{J}'_{\ti{\Omega,D}}(W)=
-\Delta^{\bot}W^{\bot}-\bar{R}(W^{\bot})
-\tilde{B}(W^{\bot})
+m\|H\|\,C_{\ti{\Omega}}(W^{\bot})-\Psi_{\ti{\Omega,D}}(W^{\bot})\nu$$ 
is the $\Omega$-Jacobi operator
 acting on sections $W$ of $T\bar{M}$ along $M$. This 
is the usual Jacobi operator
 with an extra first-order differential operator $C_{\ti{\Omega}}$ depending
 on $\Omega$ and $\bar{\nabla} \Omega$, and $\Psi_{\ti{\Omega,D}}$ a suitable
linear function. 
We define a class of vector fields of $\bar{M}$ along $\phi:D\to \bar{M}$,
\begin{equation}
\bar{\mathcal{F}}_{D,\Omega}=\{W= f(\nu+ N):
f\in\mathcal{F}_D, ~N\in C^{\infty}(\phi^{-1}T\bar{M}), ~N\bot \nu\},
\label{deformacao}
\end{equation}
and $\mathcal{F}_{D,\Omega}=\bar{\mathcal{F}}_{D,\Omega}\cap 
C^{\infty}(NM_{/D})$.
An element $W\in \bar{\mathcal{F}}_{D,\Omega}$ satisfies $W=0$ on 
$\partial D$ and 
$$\int_D\Omega(W,d\phi(e_1), \ldots, d\phi(e_m))dM=\int_D \bar{g}(W,\nu)dM=
0,$$
where $e_i$ is a direct o.n.\ frame of $M$. Such vector fields are
vectors of variation for some $\Omega$-volume preserving variations.
We will say that a submanifold with parallel mean curvature is $\Omega$-stable
if $I_{\Omega}(W,W)\geq 0$ for all  vector variations 
$W$ lying in  $ H^1_{0,T}(NM_{/D})$, i.e.,  the 
$H^1$-completion  of the vector space
generated by
$\mathcal{F}_{D,\Omega}$.
If a calibrated extension $M'$ of $M$ exists,
our stability condition is more restrictive than the one in
\cite{BdCE}, and depends on the geometry of $\bar{M}$.
But the two concepts are related, if, for example, $\Omega$ is
a parallel calibration and  $NM$ is a trivial bundle, 
or it is defined by a
fibration of $\bar{M}$ by totally geodesic
$(m+1)$-dimensional submanifolds. Related to this last case,
we  study the $\Omega$-stability of $m$-dimensional
geodesic spheres of $(m+n)$-dimensional space forms.
For the case $\bar{M}=\mathbb{R}^{m+n}$ with any parallel calibration, we
give some natural conditions in Theorem 4.2, 
which extend the case $n=1$ of \cite{BdC}
and   enable us to conclude
that a $m$-dimensional stable closed submanifold
 must be pseudo-umbilical or
even  a Euclidean sphere. 
A first difficulty in the general 
case $n\geq 2$ 
arises from the fact that a calibrated submanifold $M'$ does not have
to be totally geodesic, and stability, with no further assumptions,
 does not seem to imply this.
The Hodge theory of spheres yields other conditions on $\Omega$ that are 
necessary for their $\Omega$-stability in Euclidean spaces (Proposition 4.5).
\section{Critical area  under volume constraints}
We consider $\bar{M}$ with Riemannian metric
$\bar{g}$   and a fixed  $(m+1)$-form
 $\Omega$, and  $\phi:M\to \bar{M}$  an immersed  oriented
submanifold.
We use $\nabla$, $\nabla^{\bot}$ and $\bar{\nabla}$ to denote
the connections on $M$, $NM$ and $\bar{M}$, respectively, and
$B$ the second fundamental form of $\phi$, as a tensor with
values on  the normal bundle $NM$.

Let $D\subset M$ be a compact domain with smooth boundary, and
 $\bar{\phi}:(-\epsilon,\epsilon)\times D \to \bar{M}$,
  $\bar{\phi}(t,p)=\phi_t(p)$, a smooth variation
of $\phi=\phi_0:D\to \bar{M}$
that fixes the boundary $\forall t$. Then
the vector variation $W_{t}(p)=\frac{\partial \bar{\phi}}{\partial t}(t,p)$
vanishes at $\partial D$. If $M$ is closed ( that is, compact without boundary)
 we may consider $D=M$.  We denote by $dM_t$ the volume element of
$M_t=(M, g_t)$, where $g_t=\phi_t^*\bar{g}$, and by  $B_t$ the second
fundamental form of $M_t$.  The mean curvature vector is
$H_t =\frac{1}{m}\,  \mbox{trace}_{g_t}B_t=\sum_{ij}
\frac{1}{m}g_t^{ij}(\bar{\nabla}_{e_i}e_j)^{\bot},$ 
where $e_i$ is an oriented $g$-o.n. frame ($g=g_0$, $M=M_0$, etc..), 
$(g_t)_{ij}=g_t(e_i,e_j)$, 
and $\top$ and $\bot$ represent the orthogonal projection of 
$T\bar{M}$ onto
$TM_t$ and onto the normal bundle $NM_t$ of $M_t$, respectively.
 $M$  has parallel (constant, resp.) mean curvature 
if $H$ is a parallel section in 
the normal bundle ($\|H\|$ is constant, resp.).
The area of $D_t=(D, g_t)$ and the $\Omega$-volume of $\bar{\phi}$
are respectively given by
$$A_D(t)=\int_D dM_t,\quad \quad
V_D(t):=\int_{[0,t]\times D}\bar{\phi}^*\Omega.$$
\begin{lemma} For a local direct $g$-o.n.\ frame $e_i$ of $M$,
\begin{eqnarray*}
A_D'(t)&=&-\int_D m\bar{g}(H_t,W_t)dM_t\\
V'_D(t)&=&\int_D\Omega(W_t(p),
d\phi_t(p)(e_1), \ldots,d\phi_t(p)(e_m))dM.\end{eqnarray*}
In particular,
$A'_D(t)$ and $V'_D(t)$ depend only on $W_t^{\bot}\in NM_t$.
\end{lemma}
\noindent
\em Proof. \em
The formula for $A'_D(t)$ is very well known, but we recall here
some formulas that we will need to use  in  section 3.
 Set $g_{ij}(t,p)=g_t(e_i,e_j)$.  
Then $A_D(t)=\int_DdM_t(e_1, \ldots, e_m)dM
$ $=\int_D \sqrt{det[g_{ij}(t,p)]}dM.$ 
Using the Hessian of $\bar{\phi}$,
as a map from $ (-\epsilon, \epsilon)\times D$ with metric $ dt^2+g$,
we see that,  at $t=0$ and $p\in D$,
 $\bar{\nabla}_{\frac{d}{dt}}
(d\phi_t(e_i))= \bar{\nabla}_{e_i}W. $
Thus,
\begin{equation}
\sm{\frac{d}{dt}_{|_{t=0}}}g_{ij}=\bar{g}(\bar{\nabla}_{e_i}W, d\phi(e_j))+
\bar{g}(\bar{\nabla}_{e_j}W, d\phi(e_i)), \label{dtgij}
\end{equation}
and so 
$
\sm{\frac{d}{dt}_{|_{t=0}}}\!({det[g_{ij}(t,p)]})^{\frac{1}{2}}
=
\sum_{i}\bar{g}(\bar{\nabla}_{e_i}W, d\phi(e_i))
=
\mathrm{div}_M(W^{\top}\!)-m\bar{g}(H, W^{\bot}).
$
Therefore, ~$A'(0)=-\int_Dm\bar{g}(H, W^{\bot})dM$.
The same formula holds for any $t$. Now 
for any $~0\leq s\leq t$ and $p\in D$, 
$~\bar{\phi}^*\Omega(s,p)=$  
$\Omega( \frac{\partial \bar{\phi}}{\partial t}(s,p),
d\phi_s(p)(e_1),\ldots, d\phi_s(p)(e_m))ds\wedge dM.$ Hence
\begin{eqnarray*}
V_D(t) &=& 
 \int_0^t\La{(}\int_D\Omega( \frac{\partial \bar{\phi}}{\partial t}(s,p),
d\phi_s(p)(e_1), \ldots,d\phi_s(p)(e_m))dM\La{)} ds.
\end{eqnarray*}
Differentiation with respect to $t$ proves the lemma.\qed
\begin{definition} 
  A variation $\phi_t$  is said 
$\Omega$-volume preserving if ~$V_D(t)=V(0)=0$~  $\forall t$.
\end{definition}
For each  $W \in T_{\phi(p)}\bar{M}$  we set
\begin{equation}\label{aW}
a_W(p)=\Omega(W_p,d\phi(e_1), \ldots, d\phi(e_m))=a_{W^{\bot}}(p).\\[1mm]
\end{equation}

In what follows,  $\Omega$  is a rank-$(m+1)$
pre-calibration on $\bar{M}$.\\[2mm]
This means  $\Omega$ is an
 $(m+1)$-form  on $\bar{M}$ such that
$|\Omega(u_1, \ldots, u_{m+1})|\leq 1$,
for any o.n.\ system $u_i$ of $T_x\bar{M}$,
and equality holds for some system in $T_x\bar{M}$, at each $x\in \bar{M}$.
In the latter case, we will refer to the subspace $\mathrm{span}\{u_i\}$  
as $\Omega$-calibrated. An $\Omega$-calibrated submanifold
is an $(m+1)$-dimensional submanifold $M'$ with calibrated tangent space
(\cite{HL}). 
For these submanifolds, $\Omega$ restricted to 
$M'$ is the volume element of $M'$. If $\Omega$ is a calibration, that is, 
$\Omega$ is a closed form,
such submanifolds are homologically area minimizing,
and in particular  minimal stable in $\bar{M}$. Before we give the next
definition, we
recall the following Lemma 2.1 of
\cite{LS}:
\begin{lemma} If $u_i$ is an o.n.
 system with $\Omega(u_1, \ldots, u_{m+1})$ $=\cos\theta$, 
$\theta\in [0, \pi]$, then
for any $w\bot u_i$ $\forall i$, and any $j$,  
$|\Omega(w,u_1, \ldots, \hat{u}_j,
\ldots, u_{m+1})|\leq \sin\theta$. \end{lemma}
\noindent
\begin{definition} We will say that an oriented $m$-dimensional immersed 
submanifold $\phi:M\to \bar{M}$ has $\Omega$-calibrated extended tangent
space on $D$, if for each $p\in D$, there exist a unit normal 
vector $\nu_p$
 such that  $\bar{g}(H,\nu)=\|H\|$, 
and for a
direct o.n. frame $e_i$ of $T_pM$,
$$a_{\nu}(p)=\Omega(\nu_p,d\phi(e_1), \ldots,d\phi(e_m))=1.$$
We will refer to $EM_p=\mathbb{R}\, \nu_p\oplus T_pM$ as the extended tangent 
space
of $M$ at $p$ in $\bar{M}$, and  denote  $B^{\nu}(X,Y)=\bar{g}(B(X,Y),\nu)$.
\end{definition}
\noindent
The next lemma ensures that, if  $\nu\in NM_p$ 
 satisfies $a_{\nu}(p)=1$, then $\nu$ is unique. We will always assume
that $\nu$ defines a smooth global section of $NM_{/D}$.
\begin{lemma} If $\phi$ has calibrated extended tangent space,
then   $\forall W \in NM_p$,  $u_i\in T_pM$, 
$$\begin{array}{l}
\Omega(W,\nu, d\phi(u_1), \ldots, d\phi(u_{m-1}))=0\\
\Omega(W, d\phi(u_1), \ldots, d\phi(u_m))=\bar{g}(W, \nu_p)
\sqrt{\mathrm{det}[g(u_i,u_j)]}
\end{array}$$
Furthermore, for any $\bar{W}\in T_{\phi(p)}\bar{M}$, ~
$\bar{g}(H,\bar{W})=\|H\|a_{\bar{W}}$.
\end{lemma}
\noindent
\em Proof. \em The  equalities are immediate consequences of
Lemma 2.2. Then it follows $\bar{g}(\bar{W},H)=
\|H\|\bar{g}(\bar{W},\nu)=\|H\|a_{\bar{W}}$.\qed\\[3mm]
It is clear that if an $(m+1)$-dimensional submanifold  $M'$ of $\bar{M}$
contains $M$,
 the mean curvature of $M$
in $M'$ is the same as in $\bar{M}$  only if  $TM'=EM$ along $M$.
This is the case when $M'$ is totally geodesic
in $\bar{M}$.
We do not know if a calibrated $(m+1)$-dimensional 
submanifold $M'$ containing $M$  does exist.
Harvey and Lawson \cite{HL}, using methods of
Cartan-K\"{a}hler theory, proved that, for  some calibrations, the boundaries
of $\Omega$-calibrated manifolds are exactly the $m$-dimensional
submanifolds $\Gamma$ that are maximally $\Omega$-like, that is, 
at each $x\in \Gamma$, its
tangent space is in the span of a calibrated subspace $E_x$. 
Definition 2.2 is a particular case of
this condition. A positive 
 answer to
this  problem for a given $\Omega$ 
 would be equivalent, at each $x\in \bar{M}$, 
  to prove a modified version of  Hilbert's
 seventeenth problem in $\mathbb{R}^{m+n}\equiv T_x\bar{M}$, 
in the terms formulated in \cite{HL}.
Thus, if this problem turns out to be true for  $\Omega$, 
an $m$-dimensional  submanifold $M$ with calibrated extended tangent 
space only exists
if the mean curvature $H$ points in the same
direction of the  unit normal  of $M$ as a submanifold
of the extended calibrated manifold $M'$. 
\begin{example}\rm{ If  $\bar{M}={N}^{m+1}\times P^{n-1}$,
where $N$ and $P$ are Riemannian manifolds,
and $\Omega=\mathrm{Vol}_{N}$,
the calibrated submanifolds are the slices $N\times {h_0}$
where $h_0\in P$. Let
 $\phi:M \to \bar{M}$
with components $\phi(p)=(\psi(p), h(p))$. Then $a_{\nu}(p)=1$
for all $p$ means that $d\phi(e_i)$ and $ \nu$ lie in
$TN$. In particular $dh\equiv 0$, that is,  $h$ is constant.
Consequently $\phi$ lies in a slice
$M' =N \times h_0$. Furthermore, since $M'$ is totally geodesic
in $\bar{M}$, then $TM'=EM$ along $M$.}
\end{example}
\begin{example}  \rm{ Consider $\mathbb{R}^8$ with its octonionic structure
and $\Omega$ the Cayley calibration, $\Omega(z,u,v,w)=\bar{g}(z,u\times
v\times w)$, using the cross product of 3 vectors in $\mathbb{R}^8$
(see chapter IV  of \cite{HL}).
Let $\phi:M^3\to \mathbb{R}^8$ be any
 embedded real-analytic 3-dimensional submanifold.
Then there exists a unique 4-dimensional  real-analytic Cayley submanifold
$N$ that contains $M$ (\cite{HL}, Theorem 4.3 of Chapter IV). 
$N$ is characterized as the submanifold whose
tangent space  at $p\in N$ satisfies
 $T_pN=T_pM\oplus \mathbb{R}\mu_p$, with
 $\mu_p=d\phi(e_1)\times d\phi(e_2)\times d\phi(e_3)$.
Thus $a_{\mu}=1$, and so, 
if $M$ has 
$\Omega$-calibrated extended tangent space, then $\nu=\mu$, and $M$
has a calibrated  extension $M'=N$ such that $TM'=EM$. 
 There are many Cayley submanifolds. They can be seen as 
 the class of minimal 4-submanifolds of $\mathbb{C}^4$
 with equal K\"{a}hler angles, which includes the complex 
and the special Lagrangian  submanifolds.}
\end{example}
\begin{example} \cite{Gu1,Gu2}
\rm{Consider the imaginary part of the octonionic space
 $\mathbb{R}^7=\mathrm{Im}\,\mathbb{R}^8$, 
with o.n.\ basis $\epsilon_1, \ldots, \epsilon_7$ orthogonal
to the scalars $\mathbb{R}1$, endowed with the associative calibration
$\Omega(z,u,v)=\bar{g}(z,u\cdot v)$, where $\cdot$ is the Cayley multiplication
(see \cite{HL}, Chap.\ IV).  
Let $H:\mathbb{R}^7\times \mathbb{R}^7\to \mathbb{R}^7$ 
be the cross product of two octonions,
 $H(u,v)=\mathrm{Im}(u\cdot v)=u\times v$. If $r\leq 1$ is fixed and
$\Gamma=\{(r\cos\theta,r\sin\theta,\phi(\theta),0,0,0,0):\theta\in
[0,2\pi]\}$ is the graph of a smooth function $\phi$ over a circle
$\Gamma_0$ of radius $r$ in the $\{\epsilon_1,\epsilon_2\}$-plane, 
and if $\Gamma$ lies in a ball of radius 1,
then $\Gamma$ bounds  a surface $D$ of prescribed mean curvature 
$H(e_1,e_2)$, where $e_1,e_2$ is an o.n.\ frame of $D$ (\cite{Gu2} 3.6). 
Furthermore, $\|H\|=1$ is constant. In this case
$H=\nu$ in our setting, since $\Omega(H,e_1,e_2)=1$. Thus, $D$ is a surface
of constant mean curvature in $\mathbb{R}^7$ and  with
$\Omega$-calibrated  extended space. Moreover, if $M=D$ is
a real analytic 2-dimensional surface of
$\mathbb{R}^7$, then   $EM=TM'$, where $M'$ is the unique real analytic 
associative submanifold of $\mathbb{R}^7$ which contains $M$ (see 
Theorem 4.1 of \cite{HL}).}\end{example}
\begin{example}\rm{For any
 immersed submanifold $\phi:M^m\to\bar{M}$, one has
$\|B\|^2\geq m\|H\|^2$, and equality holds
 if and only if $\phi$ is totally umbilical (a proof can be found
in \cite{BdC} for the case $n=1$). 
This is equivalent to the second fundamental form being
 an $NM$-valued multiple
of the metric, that is,  $B(\cdot, \cdot)=H\otimes g$. 
Totally umbilical submanifolds
of space forms are also space forms and are either totally geodesic
or  $m$-spheres (see \cite{Chen}).
 If $M$ is a closed submanifold of $\mathbb{R}^{m+n}$ and 
with parallel mean curvature,  it is sufficient
to assume  $\|B\|^2\leq \frac{m^2}{m-1}\|H\|^2$ in order to conclude that 
$M$ is an $m$-dimensional sphere (\cite{ChNo}).  
A weaker concept is  pseudo-umbilicity, that is,
 when  $\bar{g}(B(\cdot, \cdot),H)=\|H\|^2g$. Chen and Yano in \cite{CY}
proved that $M$ is pseudo-umbilical with parallel mean curvature 
$H=\|H\|\nu\neq 0$
if and only if $\phi+\|H\|^{-1}\nu$ is a constant vector $z$ in 
$\mathbb{R}^{n+m}$. In this case $M$ is immersed into a hypersphere
of $\mathbb{R}^{n+m}$
centered at $z$, and $\nu$ is parallel to the radius vector field
$\phi-z$. 
We now suppose $\bar{M}=\mathbb{R}^{m+n}$ with a pre-calibration $\Omega$,
and $M$ is a  submanifold with 
 nonzero parallel mean curvature $H$ and calibrated extended
space. Let $\mathbb{R}^*=\mathbb{R}\backslash \{-\|H\|^{-1}\}$.
\begin{proposition}  If 
$M$ is pseudo-umbilical, then $\Phi: M\times \mathbb{R}^*\to \mathbb{R}^{m+n}$, 
given by $\Phi(p,t)=\phi(p)-t\nu_p$, 
 defines an  $\Omega$-calibrated  extension  $M'$ of $M$,
with second fundamental form satisfying
$B^{M'}(X',\nu)=0$ for all $X'\in T_{(p,0)}M'$, $p\in M$, and $M$
is totally umbilical in $M'$. Furthermore, supposing $M$ is closed,
then $M$ is totally umbilical in $\bar{M}$ if and only if
$M'$ is an $(m+1)$-dimensional  vector subspace and 
$M$ is a Euclidean sphere.
\end{proposition}
\noindent \em Proof. \em 
Using the pseudo-umbilicity assumption, we have for $X\in T_pM$
$$d\nu(p)(X)=\nabla^{\bot}_X\nu+\sum_i\bar{g}(d\nu(p)(X),e_i) e_i
=-B^{\nu}(X,e_i)e_i=-\|H\|X.$$
The induced metric in $M\times \mathbb{R}^*$ is
$g'_{(p,t)}=(1+t\|H\|)^2g_p+ dt^2$, and  $\Phi^*\Omega$ takes
the value $1$ along the $g'$-o.n. frame 
$\{-\frac{d}{dt},  e(t)_i=(1+t\|H\|)^{-1}e_i, i=1, \ldots, m\}$. 
Thus, $M'$ is a calibrated extension of $M$.
The tangent and the normal bundles of $M'$ at $(p,t)$ are
naturally  identified with the corresponding ones at $(p,0)$.
The global section of $TM'$, 
$\tilde{\nu}(p,t)=\nu_p$,  extends
the parallel section $\nu$ of the normal bundle of $M$,
and satisfies $d\tilde{\nu}(\tilde{\nu})=d\tilde{\nu}(-\frac{d}{dt})=0$.
Thus, the second
fundamental form  of $M'$ satisfies at $(p,t)$~ 
$B^{M'}(e_i, \tilde{\nu})=0$, and 
$B^{M'}(\tilde{\nu},\tilde{\nu})=proj_{NM'}(d\tilde{\nu}(\tilde{\nu}))=0$.
Since $B(e_i,e_j)=B^{\nu}(e_i,e_j)\nu$ $+B^{M'}(e_i,e_j)$, we conclude
that $M$ is totally
umbilical in $\bar{M}$ if and only if $B^{M'}(e_i,e_j)=0$, for all $ij$.
This holds if and only if $M'$ is totally geodesic. In this case,  $M$ is
an umbilical hypersurface of a Euclidean space, and supposing $M$ is closed, 
then by a classical result due to E.\ Cartan (or using \cite{Chen, ChNo}),
$M$ must be a sphere.
\qed
}
\end{example}
Henceforth we assume $\phi:M\to \bar{M}$ has calibrated extended
tangent space.\\[1mm]
We consider 
 the following class of functions,  defined  in \cite{BdC,BdCE}, and a class
of vector fields
$$\begin{array}{c}
\mathcal{F}_D=\{f:D\to \mathbb{R}: ~f_{/\partial D}=0,~~\int_DfdM=0\}\\[1mm]
{\mathcal{F}}_{D,\Omega}=\bar{\mathcal{F}}_{D,\Omega}
\cap C^{\infty}(NM_{/D}),\end{array}$$
where $\bar{\mathcal{F}}_{D,\Omega}$ is defined in (\ref{deformacao}).
We consider the orthogonal split of the normal bundle
$NM_{/D}=\mathbb{R}\nu\oplus F$. For each section $W\in NM$ 
we denote the corresponding
split
$$W=W^{\nu}+W^F=f\nu + W^F.$$
 \begin{definition} 
A variation $\phi_t$ of $\phi$ is said to be in $\bar{\mathcal{F}}_{D,\Omega}$ 
(resp.\ ${\mathcal{F}}_{D,\Omega}$)~ if the vector variation at $t=0$, 
$W=\frac{ \partial \bar{\phi}}{\partial t}{|_{t=0}}$, lies in 
$ \bar{\mathcal{F}}_{D,\Omega}$
(resp.\ ${\mathcal{F}}_{D,\Omega}$).
\end{definition}
\noindent
\begin{lemma} For any  
 $W\in \bar{\mathcal{F}}_{D,\Omega}$ 
 there exists an $\Omega$-volume preserving variation ${\phi}_t$ of $\phi$
that fixes the boundary and has vector variation $W$.
Reciprocally, the vector variation $W$ of any $\Omega$-volume preserving 
variation $\phi_t$ satisfies $\int_Da_WdM=0$ (not necessarily 
in $\bar{\mathcal{F}}_{D,\Omega}$).
\end{lemma}
\noindent
\em Proof. \em 
We follow the argument of \cite{BdCE}. As in (\ref{deformacao}), 
 $W=fN'$, where $N'=\nu +N$.
Let $\rho: (-\epsilon,\epsilon)\times M\to \bar{M}$ be
a  variation $\rho(\xi, p)$ such that 
$\rho(0,p)=\phi(p)$, and  $\frac{d\rho}{d\xi}(0, p)=N'_p$.
For example, we may take
$\rho(\xi,p)=exp_{\phi(p)}(\xi N'_p)$, 
where $exp$ is the exponential of $\bar{M}$. 
We consider, for each $p\in D$, the solution $\xi(t,p)$ of the
initial value problem
$$ \left\{\begin{array}{l}
\frac{d\xi}{dt}(t,p)= \frac{a_W(p)}{a(\xi(t,p),p)}\\
\ft{\xi(0,p)=0,}\end{array}\right.$$
with~ $
a(\xi,p)=\Omega(\frac{d\rho}{d\xi}(\xi,p), d\rho_{\xi}(p)(e_1),
\ldots, d\rho_{\xi}(p)(e_m))=
\sqrt{det[g_{\xi}(e_i,e_j)]} a_{\frac{\partial\rho}{\partial\xi}}$,
where $e_1, \ldots, e_m$ is any direct $g$-o.n. basis of $T_pM$.
Note that $a(0,p)=\bar{g}(N',\nu)>0$, and so for $t$ sufficiently
small $a(\xi(t,p),p)$ does not vanish.
Now  ${\phi}_t(p)=\rho(\xi(t, p),p)$ satisfies the conditions
of the lemma.
Reciprocally, if $\phi_t$ is $\Omega$-volume preserving, then 
$V'_D(0)=0$, which implies 
$0=\int_D\Omega(W_p,d\phi(e_1), \ldots, d\phi(e_m))dM=
\int_D a_WdM=\int_Df dM.$\qed\\[4mm]
The first part (a) of the next lemma is due to \cite{BdC}.
 If $D=M$ is a closed manifold,
we show  a similar conclusion for the Sobolev space
$H^1(D)=\{f\in L^2(D): \exists \nabla f\in L^2(D) 
$ ~(\mbox{in~the~weak~sense})$\}$, with the $H^1$-inner product 
$$\langle f,f'\rangle_{H^1}=\langle f,f'\rangle_{L^2} 
+\langle \nabla f,\nabla f'\rangle_{L^2}=\int_D ff' dM+\int_D
\bar{g}( \nabla f,\nabla f')dM.$$
The $L^2$-completion of $\mathcal{F}_D$ is the space $L^2_T(D)$ 
 of $L^2(D)$-functions
with zero mean value. The $H^1$-completion of $\mathcal{F}_D$ 
is $H^1_{0,T}(D)=H^1_0(D)\cap L^2_T(D)$
(see \cite{BaBe} and recall that the set of functions 
$f\in C^{\infty}(\bar{D})$
with $f_{/\partial D}=0$,  
and $\mathcal{D}(D)$ of the ones with compact support inside
$\mathring{D}$,
generate the same spaces $L^2(D)$ and  $H^1_0(D)$). If $D=M$ is closed,
$H^1_0(M)=H^1(M)$.
\begin{lemma} (a)\cite{BdC} If $G\in C^{\infty}(D)$ is $L^2$-orthogonal
to $ \mathcal{F}_D$, then $G$ is 
constant.\\[1mm]
(b) If $D=M$ is closed and   $G\in H^1(M)$ is $H^1$-orthogonal to
$\mathcal{F}_M$, then $G$ is constant a.e..
\end{lemma}
\noindent
\em Proof. \em  (b)
Let $G_M=|M|^{-1}\int_MGdM$,
where $|M|=\int_MdM$. 
Then, as $M$ is bounded, we have
$L^2(M)\subset L^1(M)$ and $G-G_M\in H^1_{0,T}(M)$.
 From $\langle G, G-G_M\rangle_{H^1}=0$, we have
$$\int_MG^2dM-G_M\int_MGdM=- \int_M\|\nabla G\|^2dM\leq 0.$$
Thus, $\int_MG^2dM\leq |M|^{-1}(\int_MGdM)^2\leq \int_MG^2dM$, where
we have used Cauchy-Schwarz in the last inequality.
Hence
 $|\langle G,1\rangle_{L^2}|$ $=|G|_{L^2}|1|_{L^2}$, which implies
$G$ is constant a.e..\qed \\[4mm]
We consider the set
\begin{equation}\label{flinha}
\mathcal{F}'_{D}= \mathcal{F}_{D}\cdot C^{\infty}(D)=
\{fh:~ f\in \mathcal{F}_{D},~ h\in C^{\infty}(D)\}
\end{equation}
spanning a vector space $\mathbb{R}\mathcal{F}'_{D}$
of the finite sums $\sum_i f_ih_i$, where
$f_i\in\mathcal{F}_{D}$ and $h_i\in C^{\infty}(D)$.
\begin{lemma}  If  $G\in L^2(D)$ (resp.\ $G\in H^1_0(D)$) 
is $L^2$-orthogonal  (resp.\ $H^1$-orthogonal) to $\mathcal{F}'_{D}$, 
then $G=0$ a.e..
\end{lemma}
\noindent
\em Proof. \em 
Since $\mathcal{D}(D)$ is $L^2$-dense in $L^2(D)$, and 
$H^1$-dense in $H^1_0(D)$, if we prove
that $\mathcal{D}(D)\subset \mathcal{F}'_D$, then we prove the lemma. 
Let $\varphi\in \mathcal{D}(D)$. We take $D'$ a domain such that 
$\mathrm{supp}\, \varphi
\subset D'\subset \bar{D}'\subset D$, and $\phi\in \mathcal{D}(D')$,
$\phi\geq 0$, and  such that
$\phi=1$ on $\mathrm{supp}\, \varphi$. Let $\phi_{\epsilon}\in \mathcal{D}(D)$
not identically zero, 
$\phi_{\epsilon}\geq 0$, and  with
compact support inside a small ball $B_{\epsilon}$ with
$\bar{B}_{\epsilon}\subset D\backslash \bar{D}'$. 
Then we have a function
$f\in \mathcal{D}(D)$ given by $\phi$ on $D'$, and 
by $-c\phi_{\epsilon}$ on
$B_{\epsilon}$, and zero away  from these sets, where $c> 0$ is the constant 
defined by $\int_D\phi dM=c\int_D\phi_{\epsilon}DM$.
Then $f\in \mathcal{F}_D$, and $\varphi=f\varphi\in 
\mathcal{F}'_D$. \qed\\[3mm]
Let $h_{D}=\frac{1}{|D|}\int_D \|H\|dM $ be the mean value of $\|H\|$.
For a variation $\phi_t$ fixing $\partial D$, define
$$
 J_{D}(t)=A_D(t) + mh_{D}V_D(t).
$$
Then $J'_D(0)=\int_Dm(-\bar{g}(H,W^{\bot})+h_Da_{W^{\bot}})dM$.
\begin{theorem} 
Consider the following statements:\\[2mm]
$(1)$ $\|H\|=h_D$, is constant in $D$.\\[1mm]
$(2)$ $A'_D(0)=0$ for all  $\Omega$-volume preserving variations
on $D$ that fix the boundary $\partial D$.\\[1mm]
$(3)$ $J'_{D}(0)=0$ for all variations on $D$ fixing the boundary 
$\partial D$.\\[1mm]
$(2')$ the same as $(2)$, and $(3')$ the same as $(3)$,
 but for $\bar{\mathcal{F}}_{D,\Omega}$
(or $\mathcal{F}_{D,\Omega}$) variations.\\[2mm]  
The statements are all equivalent.
\end{theorem}
\noindent
\em Proof. \em  
From Lemma 2.3, (1) is equivalent to $\bar{g}(H,W)=h_Da_W$
for all $W_p\in T_p\bar{M}$, 
and $p\in D$.
Now we prove $(1)\Rightarrow (2) (\Rightarrow (2'))$. 
For an $\Omega$-volume preserving variation 
with vector variation  $W$,  by  Lemma 2.1 
$A'_D(0)=\int_D-m\bar{g}(H,W)dM=-mh_{D}\int_D a_WdM$.
The latter is zero by Lemma 2.4. 
Next we prove $(2')\Rightarrow (1)$. If $W\in \bar{\mathcal{F}}_{D,\Omega}$
(or $W\in \mathcal{F}_{D,\Omega}$), we we may take $\phi_t$  a $\Omega$-volume
preserving variation with vector variation $W$ (see Lemma 2.4). 
Since $\bar{g}(W,H)=\|H\|a_W$ (by Lemma 2.3) and by
assumption $\int_D\bar{g}(H,W) dM=0$, we have 
$\int_D\|H\|a_{W}dM=0$.
Considering any function $f\in \mathcal{F}_D$ 
and $W=f\nu$, we have $a_W=f$ and  conclude that
$\int_D f\|H\|dM=0.$ Lemma 2.5(a) gives $\|H\|$  constant.
Both  $(3)\!\Rightarrow \!(2)$, $(3')\!\Rightarrow\! (2')$
 and  $(1)\!\Rightarrow\! (3)$, $(1)\!\Rightarrow\! (3')$ are obvious using
Lemmas 2.1 and 2.3.\qed
\section{The second variation}
Let $\phi:D\subset M\to \bar{M}$ be an immersion with constant  
mean curvature $H$ and 
with calibrated extended tangent space. In Theorem 2.1 we have shown 
that $\phi$ 
is a critical point of $J_D(t)$, for all variations
fixing the boundary $\partial D$.
The Laplacian for sections
in the normal bundle is given by
 $\Delta^{\bot}W^{\bot}=\sum_i\nabla^{\bot}_{e_i}
\nabla^{\bot}_{e_i}W^{\bot}-\nabla_{\nabla_{e_i}e_i}W^{\bot}$.
We use the curvature sign of $\bar{M}$,
$\bar{R}(X,Y)=-[\bar{\nabla}_X,\bar{\nabla}_Y]
+\bar{\nabla}_{[X,Y]}$, and set
$$ \begin{array}{l}
\bar{R}(W^{\bot})=
\sum_i(\bar{R}(d\phi(e_i),W^{\bot})d\phi(e_i))^{\bot}, \\[1mm]
\tilde{B}(W^{\bot})=\sum_{ij}\bar{g}(W^{\bot},B(e_i,e_j))B(e_i,e_j).
\end{array}$$
We also define a differential operator, 
$C_{\ti{\Omega}}:C^{\infty}(\phi^{-1}T\bar{M})\to 
C^{\infty}(\phi^{-1}T\bar{M})$,
given by
\begin{eqnarray}
\lefteqn{\int_D\bar{g}(C_{\ti{\Omega}}(W),W' )dM=}\\
&=& \int_D\LA{(}\sum_i\sm{\frac{1}{2}}\left(\Omega(W^{\bot},e_1, \ldots,
\nabla^{\bot}_{e_i}{W'}^{\bot}, \ldots, e_m)+\Omega({W'}^{\bot},e_1, \ldots,
\nabla^{\bot}_{e_i}{W}^{\bot}, \ldots, e_m)\right) \nonumber\\
&&\!\!\!+\sm{\frac{1}{2}}\left((\bar{\nabla}_{W^{\bot}}\Omega)
({W' }^{\bot},d\phi(e_1),\ldots,d\phi(e_m))+(\bar{\nabla}_{{W'}^{\bot}}\Omega)
({W}^{\bot},d\phi(e_1),\ldots,d\phi(e_m))\right)\LA{)}dM.\nonumber
\end{eqnarray}
If we use  the identity
$\bar{\nabla}_{e_i}W'^{\bot}=\nabla_{e_i}^{\bot}W'^{\bot}-\sum_j \bar{g}(
B(e_i,e_j),W'^{\bot}) e_j$, 
and define a vector field $X_{W,W'}$  by
$g(X_{W,W'}, e_i)=\Omega(W^{\bot}, d\phi(e_1), \ldots, W'^{\bot},
\ldots, d\phi(e_m))$, with $W'^{\bot}$ in the $i$-position, then applying
Lemma 2.3 we may write the first term of this operator as
\begin{eqnarray}
\lefteqn{\sum_i\Omega(W^{\bot},d\phi(e_1), \ldots,{\nabla}_{e_i}^{\bot}
W'^{\bot},\ldots, d\phi(e_m))=}\\
&=&\!\!\!\!\!\!\!\!
\sum_i\Omega(W'^{\bot}\!,d\phi(e_1), \ldots, {\nabla}_{e_i}^{\bot}W^{\bot}\!\!,
\ldots, d\phi(e_m))
-\bar{\nabla}_{e_i}\Omega(W^{\bot}, d\phi(e_1), \ldots, W'^{\bot},
\ldots, d\phi(e_m)) \nonumber\\
&& -\sum_i\sum_{j\neq i}\Omega(W^{\bot}, d\phi(e_1), \ldots, B(e_i,e_j), \ldots,
W'^{\bot},\ldots, d\phi(e_m)) + \mathrm{div}(X_{WW'}).\nonumber 
\end{eqnarray}
Upon integration, $\mathrm{div}(X_{WW'})$ vanishes for $W'$ with compact 
support in $\mathring{D}$. Thus, $C_{\ti{\Omega}}$ is an $L^2$-self-adjoint 
first-order differential operator,  
only depends on $C^{\infty}(NM_{/D}),$ and takes values on 
$C^{\infty}(NM_{/D})$.  If $\nabla^{\bot}\nu=0$, 
and denoting by $\nabla^F$ the connection on $F$,
by Lemma 2.3 we see that
$\Omega(W^{\bot},d\phi(e_1), \ldots,{\nabla}_{e_i}^{\bot}
W'^{\bot},\ldots, d\phi(e_m))=
\Omega(W^{F},d\phi(e_1), \ldots,{\nabla}_{e_i}^{F}
W'^{F},\ldots,$ $ d\phi(e_m))$. In this case, and if  moreover
$\bar{\nabla}\Omega=0$, then $\bar{g}(C_{\ti{\Omega}}(W),\nu)=0$
holds for all $W$.
\begin{lemma}  For a  
variation $\phi_t$ that fixes the boundary and  with vector variation
$W$,
\begin{eqnarray}
{J''}_{D}(0)&=& \int_{D}\bar{g}\La{(}
-\Delta^{\bot}W^{\bot}- \bar{R}(W^{\bot}) -\tilde{B}(W^{\bot})
+m\|H\|C_{\ti{\Omega}}(W^{\bot})
~,~ 
W^{\bot}\La{)}  dM. \label{2ndvar}
\end{eqnarray}
In particular $J''_{D}(0)$ only depends on $W^{\bot}$.
\end{lemma}
\noindent
\em Proof. \em 
 Let $e_i^t$ be a $g_t$-o.n.\  with
 $e_i^0=e_i$. We have
$$\begin{array}{lcl}
a_{W_t^{\bot}} &=&\Omega(W_t^{\bot},d\phi_t(e_1^t),\ldots, d\phi_t(e_m^t))=
\frac{\Omega(W_t^{\bot},d\phi_t(e_1),\ldots, d\phi_t(e_m))}
{(\sqrt{det[g_t(e_i,e_j)]})},\\
 J'_{D}(t)
&=&  \int_D m \La{(} -\bar{g}(H_t,W_t^{\bot})
+ h_Da_{W_t^{\bot}}\La{)}
dM_t(e_1, \ldots, e_m) dM.
\end{array}$$
By Lemma 2.3, at $t=0$, 
$\bar{g}(H_t,W_t^{\bot})= h_{D}a_{W_t^{\bot}}$. Therefore,
$$
J''_{D}(0)
= \int_D\frac{d}{dt}_{|_{t=0}}\la{(}
-m\bar{g}(H_t,W_t^{\bot})
+m h_{D}a_{W_t^{\bot}}\la{)}dM. $$
Now we have
\begin{equation}\label{var1}
\frac{d}{dt}_{|_{t=0}}\bar{g}(H_t,W_t^{\bot})=
\bar{g}(\bar{\nabla}_{\frac{d}{dt}_{|_{t=0}}}H_t,W^{\bot}) + 
\bar{g}(H, \bar{\nabla}_{\frac{d}{dt}_{|_{t=0}}}W_t^{\bot}).
\end{equation}
The next formula is well known for $W^{\top}=0$ (\cite{Si}), 
but we prove here 
the  general case 
\begin{eqnarray}\label{var2}
\bar{g}( \bar{\nabla}_{\frac{d}{dt}_{|_{t=0}}}mH_t, W^{\bot}) &=&
\bar{g}\left(\Delta^{\bot}W^{\bot} +\bar{R}(W^{\bot}) +\tilde{B}(W^{\bot})
+ m\nabla^{\bot}_{W^{\top}}H~,~ W^{\bot}\right).
\end{eqnarray}
In (\ref{var2}), if $n=1$,  the  term $\nabla^{\bot}H$  vanishes, giving the
formula  in \cite{BdC,BdCE}.
 At a fixed point 
$p_0\in M$ we consider local $g$-o.n.\ frames $e_i$ such that 
$\nabla e_i(p_0)=0$, and set $g_{ij}(t,p)=g_t(e_i,e_j)$. 
Since $mH_t=\sum_{ij}g_t^{ij}B_t(e_i,e_j)$,
then at $t=0$ and $p=p_0$
$$
  \bar{\nabla}_{\frac{d}{dt}_{|_{t=0}}}mH_t
= \sum_{ij}({\frac{d}{dt}_{|_{t=0}}}g^{ij})B(e_i,e_j)+
\sum_i\bar{\nabla}_{\frac{d}{dt}_{|_{t=0}}}(B_t(e_i,e_i)).
$$
Using the symmetry of $B$, eq.\ (2),  and  $\frac{d}{dt}_{|_{t=0}}g^{ij}
=-\frac{d}{dt}_{|_{t=0}}g_{ij}$, we have
$$\sum_{ij}({\frac{d}{dt}_{|_{t=0}}}g^{ij})B(e_i,e_j)
=\sum_i-2B(e_i,\nabla_{e_i}W^{\top})
+\sum_{ij}2\bar{g}(W^{\bot}, B(e_i,e_j))B(e_i,e_j).$$
On the other hand, at $t=0$ and $p=p_0$,
\begin{eqnarray*}
\lefteqn{\sum_i\bar{\nabla}_{\frac{d}{dt}_{|_{t=0}}}
(B_t(e_i,e_i))= \sum_i
\bar{\nabla}_{\frac{d}{dt}}\left((\bar{\nabla}_{e_i}
(d\phi_t(e_i))^{\bot}\right)}\\
&=&
\sum_i
\bar{\nabla}_{\frac{d}{dt}}\La{(}\bar{\nabla}_{e_i}
(d\phi_t(e_i))-\sum_{ku}g^{ku}_t\bar{g}(
\bar{\nabla}_{e_i}(d\phi_t(e_i)), d\phi_t(e_k))d\phi_t(e_u)\La{)}\\
&=&\sum_i \bar{\nabla}_{\frac{d}{dt}}(\bar{\nabla}_{e_i}
(d\phi_t(e_i)))
 -\sum_{iku}({\frac{d}{dt}}g^{ku}_t)
g(\nabla_{e_i}e_i(p_0),e_k)d\phi(e_u)\\
&&-\sum_{i,k}\bar{g}(
\bar{\nabla}_{\frac{d}{dt}}(\bar{\nabla}_{e_i}
(d\phi_t(e_i))),d\phi(e_k))d\phi(e_k)
-\bar{g}(\bar{\nabla}_{e_i}(d\phi(e_i)),
\bar{\nabla}_{\frac{d}{dt}}(d\phi_t(e_k)))d\phi(e_k)\\
&&-\sum_{i,k}\bar{g}(\bar{\nabla}_{e_i}(d\phi(e_i)),d\phi(e_k))
\bar{\nabla}_{\frac{d}{dt}}(d\phi_t(e_k))
\end{eqnarray*}
\begin{eqnarray*}
&=&\sum_i \La{(}\bar{\nabla}_{\frac{d}{dt}}(\bar{\nabla}_{e_i}
(d\phi_t(e_i)))\La{)}^{\bot}
-\sum_k\bar{g}(mH,\bar{\nabla}_{e_k}W)d\phi(e_k)\\
&=&\sum_i \La{(}\bar{\nabla}_{e_i}
(\bar{\nabla}_{\frac{d}{dt}}
(d\phi_t(e_i))) +\bar{R}(d\phi(e_i),W)d\phi(e_i)\La{)}^{\bot}
-\sum_k\bar{g}(mH,\bar{\nabla}_{e_k}W)d\phi(e_k)\\
&=& \La{(}\sum_i\bar{\nabla}_{e_i}\bar{\nabla}_{e_i}W
+\bar{R}(d\phi(e_i),W) d\phi(e_i)\La{)}^{\bot}
-\sum_k\bar{g}(mH,\bar{\nabla}_{e_k}W)d\phi(e_k).
\end{eqnarray*}
We note that
\begin{eqnarray*}
(\bar{\nabla}_{e_i}\bar{\nabla}_{e_i} W^{\bot})^{\bot}&=&
\La{(}\bar{\nabla}_{e_i}(\sum_j\bar{g}(\bar{\nabla}_{e_i} W^{\bot},d\phi(e_j))
d\phi(e_j)
+\nabla^{\bot}_{e_i}W^{\bot})\La{)}^{\bot}\\
&=& \La{(}\sum_j\bar{\nabla}_{e_i}(-\bar{g}(W^{\bot},B(e_i,e_j))d\phi(e_j)
)\La{)}^{\bot}
+\nabla^{\bot}_{e_i}\nabla^{\bot}_{e_i}W^{\bot}\\
&=&\sum_j -\bar{g}(W^{\bot}, B(e_i,e_j)) B(e_i,e_j)
+\nabla^{\bot}_{e_i}\nabla^{\bot}_{e_i}W^{\bot},\\
(\bar{\nabla}_{e_i}\bar{\nabla}_{e_i} W^{\top})^{\bot}&=&
\La{(}\bar{\nabla}_{e_i}({\nabla}_{e_i} W^{\top}+B(e_i,W^{\top}))
\La{)}^{\bot}\\
&=&\La{(}\bar{\nabla}_{e_i}({\nabla}_{e_i} W^{\top})\La{)}^{\bot}
+\nabla_{e_i}B(e_i,W^{\top})+ B(e_i,\nabla_{e_i}W^{\top})\\
&=& 2 B(e_i,\nabla_{e_i}W^{\top})+ \nabla_{W^{\top}}B(e_i,e_i)
-(\bar{R}(d\phi(e_i),W^{\top})d\phi(e_i))^{\bot}, 
\end{eqnarray*}
where in the last equality we have used Coddazzi's equation.
Here $\nabla B$ denotes the covariant derivative of $B$ as a tensor
with values in $NM$.
Therefore,
\begin{equation}\label{var3}
\bar{\nabla}_{\frac{d}{dt}_{|_{t=0}}}mH_t=
\Delta^{\bot}W^{\bot}+\tilde{B}(W^{\bot})+\bar{R}(W^{\bot})
+\nabla_{W^{\top}}^{\bot}mH
-\sum_k\bar{g}(mH,\bar{\nabla}_{e_k}W)d\phi(e_k)
\end{equation}
and we  obtain (\ref{var2}). 
 Next we calculate $\frac{d}{dt}_{|_{t=0}}a_{W_t^{\bot}}$.
\begin{eqnarray*}
\lefteqn{\frac{d}{dt}_{|_{t=0}}a_{W_t^{\bot}}= 
\frac{d}{dt}_{|_{t=0}}\!\!\!\!\!
\Omega(W_t^{\bot},d\phi_t(e_1),\ldots, d\phi_t(e_m)) +
\bar{g}(W^{\bot}, \nu)\frac{d}{dt}_{|_{t=0}}\!\!\!\!\!(det[g_t(e_i,e_j)])
^{-\frac{1}{2}}}\\
&=&(\bar{\nabla}_{\frac{d\phi_t}{dt}_{|_{t=0}}}\Omega)(W^{\bot},d\phi(e_1),
\ldots, d\phi(e_m))
+ \Omega(\bar{\nabla}_{\frac{d}{dt}_{|_{t=0}}}\!\!\!\!\!W_t^{\bot},
 d\phi(e_1),\ldots, d\phi(e_m))\\
&&+\sum_i \Omega(W^{\bot}, d\phi(e_1), \ldots, 
\bar{\nabla}_{\frac{d}{dt}_{|_{t=0}}}\!\!\!\!\!(d\phi_t(e_i)), \ldots,
 d\phi(e_m))-\bar{g}(W^{\bot}, \nu)
\bar{g}(\bar{\nabla}_{e_i}W,d\phi(e_i)).
\end{eqnarray*}
Applying Lemma 2.3, and  
$\bar{\nabla}_{\frac{d}{dt}_{|_{t=0}}}
(d\phi_t(e_i))=\bar{\nabla}_{e_i}W= (\bar{\nabla}_{e_i}W)^{\top}
+\nabla^{\bot}_{e_i}W^{\bot}+B(e_i,W^{\top})$, 
we have
\begin{eqnarray*}
\Omega(W^{\bot}, d\phi(e_1), \ldots,
\bar{\nabla}_{\frac{d}{dt}_{|_{t=0}}}\!\!\!\!\!(d\phi_t(e_i)), \ldots,
 d\phi(e_m)) ~=~
\bar{g}(W^{\bot},\nu)\bar{g}(\bar{\nabla}_{e_i}W,d\phi(e_i))~+ \\
+ \sum_i\Omega(W^{\bot},d\phi(e_1), \ldots,\nabla^{\bot}_{e_i}W^{\bot},
\ldots, d\phi(e_m))+\Omega(W^{\bot},d\phi(e_1), 
\ldots,B(e_i,W^{\top}\!\!),\ldots,d\phi(e_m)).
\end{eqnarray*}
Thus,
\begin{eqnarray}\label{var4}
\frac{d}{dt}_{|_{t=0}}a_{W_t^{\bot}} &=& 
(\bar{\nabla}_{W^{\top}}\Omega)(W^{\bot},d\phi(e_1),\ldots, d\phi(e_m))+
\bar{g}(C_{\ti{\Omega}}(W^{\bot}),
W^{\bot})\\
&&+ \bar{g}(\bar{\nabla}_{\frac{d}{dt}_{|_{t=0}}}\!\!\!\!\!W_t^{\bot}, 
\nu) +
\sum_i\Omega(W^{\bot},d\phi(e_1), \ldots,B(e_i,W^{\top}\!),
\ldots,d\phi(e_m)). \nonumber
\end{eqnarray}
Now we observe that,
for $X\in T_{p_0}M$,
\begin{eqnarray}
\lefteqn{(\bar{\nabla}_{X}\Omega)(W^{\bot},d\phi(e_1),
\ldots, d\phi(e_m))=d(\bar{g}(W^{\bot},\nu))(X)
-\Omega(\bar{\nabla}_XW^{\bot},
d\phi(e_1),\ldots, d\phi(e_m))~~~}\nonumber\\
&&-\sum_i\Omega(W^{\bot},d\phi(e_1),
\ldots, \bar{\nabla}_X(d\phi(e_i)), \ldots, d\phi(e_m))\nonumber\\
&=& d(\bar{g}(W^{\bot},\nu))(X)-\bar{g}(\bar{\nabla}_XW^{\bot}, \nu)
-\sum_i\Omega( W^{\bot},d\phi(e_1),\ldots, B(X,e_i), 
\ldots, d\phi(e_m)).\nonumber
\end{eqnarray}
Consequently, 
\begin{eqnarray} \label{var6}
\lefteqn{(\bar{\nabla}_{X}\Omega)(W^{\bot},d\phi(e_1),
\ldots, d\phi(e_m))=}\\
&=&\bar{g}(W^{\bot},\nabla_X^{\bot}\nu)-\sum_i\Omega( W^{\bot},
d\phi(e_1),\ldots, B(X,e_i), \ldots, d\phi(e_m)).\nonumber
\end{eqnarray}
Therefore, taking $X=W^{\top}$ in (\ref{var6}), using (\ref{var4}) and 
(\ref{var2}), and adding (\ref{var1}),
\begin{eqnarray*}
\frac{d}{dt}_{|_{t=0}}\!\!\!\!\!\La{(}mg(H_t,W^{\bot}_t)
-m \|H\| a_{W_t^{\bot}}\La{)}=
 \bar{g}\la{(}\Delta^{\bot}W^{\bot}+\bar{R}(W^{\bot})
+\tilde{B}(W^{\bot})-m \|H\| C_{\ti{\Omega}}(W^{\bot})
,W^{\bot}\la{)}.\qed
\end{eqnarray*}
\noindent
From (\ref{var6})  in the preceding proof, we conclude:
\begin{proposition} If $\phi:M\to \bar{M}$ is an immersion
with calibrated extended
tangent space, then 
 $\nu$ is a parallel section of the normal bundle
if and only if, $~\forall X\in T_pM, W^{\bot}\in NM_p,$ 
$$(\bar{\nabla}_{X}\Omega)(W^{\bot},d\phi(e_1),
\ldots, d\phi(e_m))= -\sum_i\Omega( W^{\bot},
d\phi(e_1),\ldots, B(X,e_i), \ldots, d\phi(e_m)).$$
In this case, if $\phi$  has constant mean curvature, then it
has parallel mean curvature. 
\end{proposition}
We now
define  a self-adjoint strongly elliptic second order differential operator
 $\mathcal{J}_{\ti{ \Omega}}:C^{\infty}(\phi^{-1}T\bar{M})\to 
C^{\infty}(\phi^{-1}T\bar{M}),$
$$
\mathcal{J}_{\ti{\Omega}}(W)=\mathcal{J}(W)+m\|H\|\, C_{\ti{\Omega}}(W)
~~~\in C^{\infty}(NM_{/D}),
$$
where  $\mathcal{J}(W)$ is the usual Jacobi operator,
$\mathcal{J}(W)=-\Delta^{\bot}W^{\bot}-\bar{R}(W^{\bot})
-\tilde{B}(W^{\bot})$.\\

We now recall some properties of calibrations defined by fibrations.
Consider $\pi:\bar{M}\to N$  a Riemannian submersion between Riemannian
manifolds, defining an orthogonal split of $T\bar{M}$ into the vertical
and the horizontal spaces, $T\bar{M}=T\bar{M}^{\upsilon}\oplus T\bar{M}^h$.
For $y\in N$, $M'_y=\pi^{-1}(y)$ is the fibre at $y$, which we assume to
be of dimension $m+1$, and for 
$x\in M'_y$, $T_x\bar{M}^{\upsilon}=T_x(M'_y)$, $T_x\bar{M}^h=(NM'_y)_x$.
For each vector $X\in T\bar{M}$, we denote by $X^{\upsilon}$ and $X^h$
its projection into $T\bar{M}^{\upsilon}$ and $ T\bar{M}^h$, respectively.
This fibration defines a pre-calibration on $\bar{M}$ that calibrates
the fibres $M'_y$. It is given by
$$\Omega_{\pi}(X_1,\ldots,X_{m+1})=\mathrm{Vol}_y(X_1^{\upsilon},
\ldots,X_{m+1}^{\upsilon}),~~~~~\forall X_i\in T_x\bar{M}$$
where $\mathrm{Vol}_y$ is the volume element of the fibre $M'_y$, with 
$y=\pi(x)$.
Let $e'_i$, $i=1, \ldots, {m+1}$ and $e'_{\alpha}$, $\alpha=m+2,\ldots, m+n$
be local o.n.\ frames of $T\bar{M}^{\upsilon}$ and $T\bar{M}^h$, respectively.
\begin{lemma}\cite{LS} 
All components of $\bar{\nabla}\Omega_{\pi}$ and
of  $d\Omega_{\pi}$ vanish
except for the following, where $i,j\leq m+1$, $\alpha,\beta\geq m+2$:
$$ \begin{array}{l}
\bar{\nabla}_{e'_j}\Omega_{\pi} (e'_{\alpha},e'_1, \ldots, \hat{e}'_i, 
\ldots, e'_{m+1})=(-1)^{i+1}\bar{g}(B^{\upsilon}(e'_j, e'_i), e'_{\alpha}) \\
\bar{\nabla}_{e'_{\beta}}\Omega_{\pi} ( e'_{\alpha},
e'_1, \ldots, \hat{e}'_i, \ldots, e'_{m+1})
=(-1)^{i}\bar{g}(\bar{\nabla}_{e'_{\beta}}e'_{\alpha}, e'_i)\\
d\Omega_{\pi} (e'_{\alpha}, e'_1,\ldots, e'_{m+1}) =- (m+1)
\, \bar{g}(H^{\upsilon}, e'_{\alpha})\\
d\Omega_{\pi} (e'_{\alpha}, e'_{\beta},e'_1,\ldots,\hat{e}'_i,\ldots,e'_{m+1})=
(-1)^i\bar{g}([e'_{\alpha},e'_{\beta}], e'_i), 
\end{array}$$
where $B^{\upsilon}$ and $H^{\upsilon}$ denote the second 
fundamental form and the mean curvature of the fibres, respectively.
\end{lemma}
\noindent
 \begin{proposition} Assume  $M'$ is a totally geodesic
fibre of a Riemannian submersion $\pi:\bar{M}\to N$. Furthermore,
assume  $\phi:M\to \bar{M}$ is an immersion with
 $\Omega_{\pi}$-calibrated extended tangent space  and that
$\phi(M)$ lies in  $ M'$ with  $EM=TM'$ along $M$.
Then $\nu$ is a parallel section of $NM$
 and $C_{\ti{\Omega}}=0$. In particular, if $\phi$ has constant mean 
curvature, then it has parallel mean curvature. 
\end{proposition}
\noindent
\em Proof. \em We take frames $e'_a$ such that, at $p\in M$, $e'_1=\nu$,
 $e'_{i+1}=d\phi(e_i)$, for $i=1,\ldots, m$. The first equality 
of Lemma 3.2
and $B^{\upsilon}=0$ give us  $\bar{\nabla}_{e_i}\Omega_{\pi}
(e'_{\alpha},d\phi(e_1),\ldots, d\phi(e_m))=0$. 
By the Lemma, the component $\bar{\nabla}_{e_i}\Omega_{\pi}(\nu,d\phi(e_1),
\ldots, d\phi(e_m))$ also vanishes. 
It is clear that $\Omega_{\pi}(W^{\bot}, W'^{\bot}, d\phi(e_1), \ldots,
d\phi(\hat{e}_i), \ldots, d\phi(e_m))=0$. 
Applying Proposition 3.1, we  conclude
$\nu$ is  parallel in  $NM$. To prove that $C_{\ti{\Omega}}=0$ we  
use the second equality of Lemma 3.2 and the Escobales-O'Neill identity 
$~(\bar{\nabla}_{e'_{\alpha}}e'_{\beta})^{\upsilon}=
\frac{1}{2}[e'_{\alpha},e'_{\beta}]^{\upsilon}$, ~and that~
$\bar{\nabla}_{\nu}\Omega_{\pi}(e'_{\alpha},d\phi(e_1),
\ldots, d\phi(e_m))=\bar{\nabla}_{e'_a}\Omega_{\pi}(\nu,d\phi(e_1),
\ldots, d\phi(e_m))=0$, for any $a=1, \ldots, m+n$.\qed
\section{$\Omega$-stable submanifolds with parallel  mean curvature}
Let $\phi:M\to \bar{M}$ be an immersed submanifold with
calibrated extended tangent space and parallel mean curvature.
Given  a section  $W\in\bar{\mathcal{F}}_{D,\Omega}$, 
by Lemma 2.4 there is an $\Omega$-volume preserving variation of $\phi$
with vector variation $W$. For such a variation
we have $A'_D(0)=J'_D(0)=0$ and
$A''_D(0)= J''_D(0) =\int_D\bar{g}(\mathcal{J}_{\ti{\Omega}}(W^{\bot}),
W^{\bot})dM.$ We
define a symmetric bilinear operator  on 
the vector space  $\mathbb{R} \bar{\mathcal{F}}_{D,\Omega}$
spanned by $\bar{\mathcal{F}}_{D,\Omega}$
$$I_{\Omega}(W,W'):=
\int_D \bar{g}(\mathcal{J}_{\ti{\Omega}}(W),W')dM 
=I_{\Omega}(W^{\bot},W'^{\bot}).$$
We consider the orthogonal split $NM=\mathbb{R}\nu\oplus F$ into two parallel
subbundles.
 For $f\in\mathcal{F}_D$ and $W^F\in C^{\infty}(F)$, 
we have $f\nu$, $f(\nu+W^F)$ $\in \mathcal{F}_{D,\Omega}$. Then $fW^F\in
\mathbb{R}\mathcal{F}_{D,\Omega}$. Hence, 
$\mathbb{R}\mathcal{F}_{D,\Omega}
=\mathcal{F}_D\oplus \mathcal{F}'(F)$, where $f\in\mathcal{F}_D$
is identified with $f\nu$, and 
$$\mathcal{F}'(F)=\{
\sum_a{f}_aW^F_a (\mbox{finite~sum}):
{f}_a\in\mathcal{F}_D, ~W^F_a\in C^{\infty}(F)\}.$$
Let $L^2(NM_{/D})$ be the space of measurable sections $W$
of the normal bundle such that $\|W\|\in L^2(D)$, and $H^1(NM_{/D})$
the space of sections $W\in L^2(NM_{/D})$ such that
$\exists\nabla^{\bot}_X W\in L^2(NM_{/D})$ (in the weak sense)
 for all $X\in C^{\infty}(TM_{/D})$.
We  define $L^2_T(NM_{/D})$ as the $L^2$-completion of 
$ \mathbb{R}\mathcal{F}_{D,\Omega}$ in $L^2(NM_{/D})$, 
and $L'^2(F)$ the $L^2$-completion of $\mathcal{F}'(F)$.
Then $L^2_T(NM_{/D})=L^2_T(D)\oplus L'^2(F)$.
We also denote by  $H'^1_0(D)$,
 $H^1_{0,T}(NM_{/D})$, $H'^1_0(F)$ the corresponding
$H^1$-completion of  $\mathbb{R}\mathcal{F}'_D$,
 $\mathbb{R}\mathcal{F}_{D,\Omega}$, and $\mathcal{F}'(F)$, 
respectively, where
$$\langle W,W'\rangle_{H^1}=\int_D\bar{g}(W,W')dM+\int_M
\sum_i\bar{g}( \nabla_{e_i}^{\bot} W,\nabla_{e_i}^{\bot} W')
dM.$$
 If $D=M$ closed, $H^1_0=H^1$.
We consider the quadratic form defined for $W\in
 \mathbb{R}\mathcal{F}_{D,\Omega}$
\begin{equation}\label{quadraticform}
Q_{\Omega}(W)=\int_D\left(\|\nabla^{\bot}W\|^2-\bar{g}(\bar{R}(W),W)
-\bar{g}(\tilde{B}(W),W)+m \|H\|\bar{g}(C_{\Omega}(W),W)\right)dM. 
\end{equation}
Then $I_{\Omega}(W,W)=Q_{\Omega}(W)$,
and so $I_{\Omega}$ has a natural extension to 
$W\in H^1_{0,T}(NM_{/D})$.
\begin{lemma} If $Z\in L^2(NM_{/D})$ satisfies
$\int_D\bar{g}(Z,W')dM=0$~ for all $W'\in \mathbb{R} \mathcal{F}_{D,\Omega}$, 
then $~Z=c\nu$ a.e.,  where $c$ is a constant.
\end{lemma}
\noindent
\em Proof. \em If we take $W'=f\nu$ where $ f\in \mathcal{F}_D$,
then we conclude by Lemma 2.5(a) that $\bar{g}(Z,\nu)=c$ a.e.\ where
$c$ is constant. We also have $\int_Df\bar{g}(Z,W^F)dM=0$ for all
$W^F\in C^{\infty}(F)$ and $f\in \mathcal{F}_D$.
 Thus $\bar{g}(Z,W^F)$ is constant a.e.\
Taking a non-constant function $\rho$ we conclude $\bar{g}(Z,\rho W^F)$
is also constant a.e..
This implies $Z^F=0$ a.e..\qed
\begin{lemma} $C_0^{\infty}(F)=L'^2(F)\cap C_0^{\infty}(F)$.
In particular $\mathcal{F}'(F)$ is $L^2$-dense in $C_0^{\infty}(F)$.
\end{lemma}
\noindent
\em Proof. \em First we claim that if $Z\in L^2(F)$ and $Z\bot \mathcal{F}'(F)$
then $Z=0$ a.e.. To see this,  we fix $W\in L^2(F)$. Then, for any
$f\in \mathcal{F}_{D}$, $0=\int_D \bar{g}(Z,fW)dM$, which implies by 
Lemma 2.5(a)
that $\bar{g}(Z,W)$ is constant a.e.. Since $W$ is arbitrary,
  $Z=0$ a.e..Therefore, $L'^2(F)=L^2(F)$.
On the other hand, the $L^2$-closure of $C_0^{\infty}(F)$ is $L^2(F)$.
\qed
\begin{lemma}
$\mathbb{R}\mathcal{F}'_D$ is $H^1$-dense in $H^1_0(D)$, that is, 
 $H'^{1}_0(D)=H^1_0(D)$. \end{lemma}
\noindent
\em Proof. \em This is an immediate consequence of Lemma 2.6.\qed
\begin{proposition}
 $\mathcal{F}'(F)$ is $H^1$-dense in $H^1_0(F)$, that is,
 $H'^1_0(F)=H^1_0(F)$.  Furthermore,
$H^1_{0,T}(NM_{/D})=H^1_{0,T}(D)\oplus H^1_0(F)$ and, for 
$W\in H^1_0(NM_{/D})$, $W$ lies in $ H^1_{0,T}(NM_{/D})$ 
if and only if  $\int_D a_WdM=0$.
\end{proposition}
\noindent
\em Proof. \em  We only need to prove that any $W\in C_0^{\infty}(F)$ with
compact support $K\subset \mathring{D}$ is an element of $\mathcal{F}'(F)$, 
since the set
of such sections is $H^1$-dense in $H^1_0(F)$ (see \cite{Smale}).
Let $\varphi\in \mathcal{D}(D)$ with $\varphi=1$ on $K$. 
We have proved in the proof of Lemma 2.6 that $\varphi\in \mathcal{F}_D'$,
say $\varphi=fh$, as in (4).
Then $W=\varphi W=f(hW)\in \mathcal{F}'(F)$.
 The rest is elementary.\qed
\begin{remark}\rm{Assume $M$ is closed and consider $D=M$. 
In Lemma 4.3 we have shown that
the $H^1$-closure of $\mathbb{R}\mathcal{F}'_M$ is $H^1(M)$.
Let $\phi_i$, $i=0,1\ldots,$ be an $L^2$-orthonormal basis of $L^2(M)$ 
of eigenfunctions of $-\Delta$,
with corresponding eigenvalues $\lambda_i\nearrow +\infty$,
where $\lambda_0=0$ and $\lambda_1>0$. 
Then $\langle \phi_i,\phi_j\rangle_{L^2} =\delta_{ij}$ and
$\langle \nabla\phi_i,\nabla\phi_j\rangle_{L^2} =\lambda_i\delta_{ij}$.
Using integration by parts, i.e.,  $\int_M\|\nabla f\|^2=-\int_Mf\Delta f dM$,
we obtain for all $i,j$,
$$\int_M\|\nabla(\phi_i\phi_j)\|^2dM=(\lambda_i+\lambda_j)
\int_M\phi_i^2\phi_j^2 dM-2\int_M\phi_i\phi_j
{g}(\nabla \phi_i,\nabla \phi_j)dM.$$
On the other hand,
$\|\nabla(\phi_i\phi_j)\|^2=
\phi_i^2\|\nabla \phi_j\|^2+\phi_j^2\|\nabla \phi_i\|^2
+2\phi_i\phi_jg(\nabla \phi_i,\nabla \phi_j).$
Hence, 
$$\int_M\|\nabla(\phi_i\phi_j)\|^2dM=\frac{(\lambda_i+\lambda_j)}{2}
\int_M\phi_i^2\phi_j^2 dM+\frac{1}{2}\int_M(\phi_i^2\|\nabla \phi_j\|^2
+\phi_j^2\|\nabla \phi_i\|^2)dM.$$
Consequently, 
$$\int_M\|\nabla(\phi_i\phi_j)\|^2dM\geq \frac{(\lambda_i+\lambda_j)}{2}
\int_M(\phi_i\phi_j)^2dM, ~~~
\int_M\|\nabla \phi_i^2\|^2dM=\frac{4}{3}\lambda_i\int_M\phi_i^4dM.$$
If we take $i\geq 1$ and $j\geq 0$,  $\phi_i\phi_j\in\mathcal{F}'_M$
 satisfy the above inequalities. But  the constant function $h=1\in H^1(M)$
can be expressed as an $L^2$-limit of series in terms of $\phi_i\phi_j$,
and it does not satisfy an inequality 
$\int_M\|\nabla h\|^2dM\geq c\int_Mh^2dM$, where $c$ is a positive constant.
We note that $\phi_i\phi_j$, with $i\leq j$, is not an orthonormal system.
}
\end{remark}

\begin{definition} We will say an immersed submanifold $\phi:M\to \bar{M}$ of
calibrated extended tangent space and of parallel
mean curvature  is  essentially
 $\Omega$-stable on $D$ if ${A''}_D(0)\geq 0$
for all $\Omega$-volume preserving $\bar{\mathcal{F}}_{D,\Omega}$-variations. 
Equivalently, $\phi$ is essentially $\Omega$-stable on $D$ if and only if
$J''_D(0)\geq 0$ for any variation $\phi_t$ with vector variation
$W\in \bar{\mathcal{F}}_{D,\Omega}$. 
We will say $\phi$ is 
$\Omega$-stable on $D$, if  $\forall W\in H^1_{0,T}(NM_{/D})$,
$I_{\Omega}(W,W)\geq 0$,
and $\Omega$-unstable if otherwise.
\end{definition}
\noindent The equivalence of ${A''}_D(0)\geq 0$
with the  condition $J''_D(0) \geq 0$
comes from the fact that $J''_D(0)$ does not depend on the
variation $\bar{\phi}$ but only on the normal component $W^{\bot}$ of the
vector variation. The variation $\bar{\phi}$ does not need to be
$\Omega$-volume preserving, but one of the variations with
vector variation $W\in \bar{\mathcal{F}}_{D,\Omega}$
is $\Omega$-volume preserving.\\[-2mm]

We consider the linear function on $\mathbb{R}\mathcal{F}_{D,\Omega}$, 
$~\Psi_{\Omega,D}(W)=|D|^{-1}\int_D\bar{g}(\mathcal{J}_{\ti{\Omega}}(W),
\nu)dM$, and  define a self-adjoint operator, the
$\Omega$-Jacobi operator, $\mathcal{J}'_{\ti{\Omega,D}}:
\mathbb{R}\mathcal{F}_{D,\Omega}\subset H^1_{0,T}(NM_{/D}) 
\to L^2_{T}(NM_{/D})$,
given by
$$\mathcal{J}'_{\ti{\Omega,D}}(W)=\mathcal{J}_{\ti{\Omega}}(W)-\Psi_{\Omega,D}
(W)\nu.$$
 Then for all $W,W'\in \mathbb{R}\mathcal{F}_{D,\Omega}$,
$ I_{\Omega}(W,W')=\int_D\bar{g}(\mathcal{J}'_{\ti{\Omega,D}}(W),W')dM.$
We extend the
definition of Jacobi field  given in \cite{BdCE}:
\begin{definition}  We will say 
that $W\in H^1_{0,T}(NM_{/D})\cap C^{\infty}(NM_{/D})$
is  an $\Omega$-Jacobi field along $\phi:D\to \bar{M}$ if $I_{\Omega}(W,W')=0$, 
$\forall W'\in \mathbb{R}\mathcal{F}_{D,\Omega}$.
\end{definition}
 The next proposition follows immediately from the 
previous lemmas of this section:
\begin{proposition}  
$W$ is an $\Omega$-Jacobi field if and only if
$\mathcal{J}_{\ti{\Omega}}(W)=c\nu$, where $c$ is a constant,
if and only if $\mathcal{J}'_{\ti{\Omega,D}}(W)=0$.
\end{proposition}
If $\phi$ is a minimal immersion,  and $Z$ is a Killing
vector field of $\bar{M}$, it is well known that
$Z^{\bot}$ is a Jacobi field for the usual Jacobi operator 
($C_{\ti{\Omega}}=0$ ) in the sense that
$\mathcal{J}(Z^{\bot})=0$ \cite{Si}.
A proof can be obtained by recalling
that Killing vector fields generate a one-parameter family
of isometries $\Phi_t$ on $\bar{M}$, defining a variation
$\phi_t=\Phi_t\circ \phi$ by minimal immersions, and so a Jacobi field  with
vector variation.
This is also true
 if $\phi$ has
constant mean curvature  with $\bar{M}=M' $, $n=1$ and $\Omega$ is
the volume form of $M' $  \cite{BdCE}.
In higher codimension we need some additional assumptions.
\begin{proposition} If $\phi:M\to \bar{M}$ is any immersion  and $Z$ 
a Killing vector field of $\bar{M}$, then
$$
\mathcal{J}_{\ti{\Omega}}(Z^{\bot})={\nabla}_{Z^{\top}}^{\bot}mH
-(\bar{\nabla}_{mH}Z)^{\bot}+m\|H\|\, C_{\ti{\Omega}}(Z^{\bot}).
$$
Furthermore, suppose $\phi$ has extended calibrated tangent space and
a minimal calibrated extension $M'$ of $M$ exists such that
$M$ is a closed hypersurface in $M'$ as the boundary of an open domain $O'$
of $M'$.  Then $\int_{M}a_{Z^{\bot}}dM=0$.
In this case, if $\phi$ has parallel mean curvature,
then $Z^{\bot}$ is an $\Omega$-Jacobi field along $\phi$ 
 if and only if 
\begin{equation}
-(\bar{\nabla}_{m H}Z)^{\bot}+m\|H\|\, C_{\ti{\Omega}}(Z^{\bot})=c\nu,
\end{equation}
where $c= m\|H\|\bar{g}(C_{\ti{\Omega}}(Z^{\bot}),\nu)$ is a constant,
which is zero  if $\bar{\nabla}\Omega=0$.
\end{proposition}
\begin{remark}\rm{(1)
If $n=1$, then  $H$ is parallel, and 
$(\bar{\nabla}_{H}Z)^{\bot}=\|H\|
\bar{g}(\bar{\nabla}_{\nu}Z,\nu)\nu=0$,
since $Z$ is Killing. If $C_{\ti{\Omega}}=0$ 
(for example $\Omega=d\bar{M}$) then
 $\mathcal{J}_{\ti{\Omega}}(Z^{\bot})=0$.
}\end{remark}
\noindent
\em Proof. \em
 Let $p_0\in M$ and $e_i$, $W_{\alpha}$ local o.n. frames of $TM$ and $NM$,
defined on an open set $D$ of $M$ which contains $p_0$, and 
such that $\nabla_X e_i(p_0)=\nabla_X^{\bot}W_{\alpha}(p_0)=0$,
for all $X\in T_{p_0}M$. A tubular neighbourhood $\mathcal{V}$
 of $D$ in $\bar{M}$
is diffeomorphic to an open set of $NM$, using the exponential map of 
$\bar{M}$.
Each point $q\in \mathcal{V}$ is of the form $q=exp_p(v)$ for a unique
$p\in D$ and $v\in NM_p$. Let $\gamma(t)$ be the geodesic
starting at $p$ with initial velocity $v$. Then we define
$\bar{e}_i(q)$ and $\bar{W}_{\alpha}(q)$ as the parallel 
transport  along $\gamma(t)$ of $e_i(p)$ and of $W_{\alpha}(p)$, respectively.
In this way we have vector fields on $\bar{M}$
defined on a neighbourhood of $p_0$, extending $e_i$ and $W_{\alpha}$.
At $p_0$ we have
$$
 (\bar{\nabla}_{e_i}W_{\alpha})^{\top} =
-\sum_j \bar{g}(W_{\alpha},
B(e_i,e_j))e_j,~~~~ 
\bar{\nabla}_{W_{\beta}}\bar{e}_i=\bar{\nabla}_{W_{\beta}}
\bar{W}_{\alpha}=0\\[-5mm]
$$
\begin{eqnarray*}
\sum_i B(e_i, \nabla_{e_i}Z^{\top})&=&\sum_{ij}B(e_i,e_j)g(e_j,
\bar{\nabla}_{e_i}Z^{\top})\nonumber\\
&=&\sum_{ij}B(e_i,e_j)\bar{g}(d\phi(e_j),
\bar{\nabla}_{e_i}Z)-\sum_{ij}B(e_i,e_j)
\bar{g}(d\phi(e_j),\bar{\nabla}_{e_i}Z^{\bot})\nonumber\\
&=&
\sum_{ij}B(e_i,e_j)\bar{g}(B(e_i,e_j),Z^{\bot})=\tilde{B}(Z^{\bot}),\nonumber
\end{eqnarray*}
where in the last equality we have used the fact that $B(e_i,e_j)$ is
symmetric and $\bar{g}(\bar{\nabla}_{e_i}Z,e_j)$ skew-symmetric in
$ij$.  Now we have at $p_0$ (and identifying $e_i$ with $d\phi(e_i)$), 
\begin{eqnarray*}
\lefteqn{\Delta^{\bot}Z^{\bot}= \sum_{i\alpha}\left(\bar{\nabla}_{e_i}(
\bar{g}(\bar{\nabla}_{e_i}Z^{\bot}, W_{\alpha})W_{\alpha})\right)^{\bot}=
 \sum_{i\alpha}(d_{e_i}
\bar{g}(\bar{\nabla}_{e_i}Z^{\bot}, W_{\alpha}))W_{\alpha}}\\
&=& \sum_{i\alpha}\left(d_{e_i} (
\bar{g}(\bar{\nabla}_{e_i}Z, W_{\alpha})-
\bar{g}(\bar{\nabla}_{e_i}Z^{\top}, W_{\alpha}))\right)W_{\alpha}\\
&=& \sum_{i\alpha}\left(d_{e_i} (
-\bar{g}(\bar{\nabla}_{\bar{W}_{\alpha}}Z, {e}_i)-
\bar{g}(B(e_i, Z^{\top}), W_{\alpha}))\right)W_{\alpha}\\
&=& \sum_{i\alpha}\left(
-\bar{g}(\bar{\nabla}_{e_i} \bar{\nabla}_{\bar{W}_{\alpha}}Z, {e}_i)
-\bar{g}( \bar{\nabla}_{\bar{W}_{\alpha}}Z,\bar{\nabla}_{e_i}e_i)
-\bar{g}({\nabla}_{e_i}^{\bot}(B(e_i, Z^{\top})), W_{\alpha})\right)W_{\alpha}.
\end{eqnarray*}
Note that for $X,Y$ vector fields on $\bar{M}$,
 $\bar{g}(\bar{\nabla}_{X}\bar{\nabla}_{Y}Z, Y)=0$. Therefore
$
\bar{g}(\bar{\nabla}_{e_i} \bar{\nabla}_{\bar{W}_{\alpha}}Z, {e}_i)
=\bar{g}(\bar{R}(W_{\alpha}, e_i)Z +\bar{\nabla}_{[\bar{W}_{\alpha},
\bar{e}_i]}Z, e_i).$
Using  the vanishing properties of the
covariant derivatives of $e_i$ and $W_{\alpha}$ at $p_0$,
 and the fact that $B(e_i,e_j)$ is symmetric and 
$\bar{g}(\bar{\nabla}_{e_i}Z, e_j)$ is skew-symmetric in $ij$, we have 
$\sum_i\bar{g}(\bar{\nabla}_{[\bar{W}_{\alpha},
\bar{e}_i]}Z, e_i)=\sum_i-\bar{g}((\bar{\nabla}_{e_i}Z)^{\top}, 
\sum_j\bar{g}(W_{\alpha}, B(e_j,e_i))e_j)=0.$
Applying  Coddazzi's equation,
$\sum_{i}{\nabla}_{e_i}^{\bot}
(B(e_i, Z^{\top}))=\sum_i\nabla_{Z^{\top}}B(e_i,e_i)
-(\bar{R}(e_i,Z^{\bot})e_i)^{\bot}+ B(e_i,\nabla_{e_i}Z^{\top}),$
and  we  arrive at
$$\Delta^{\bot}Z^{\bot}= \sum_i -(\bar{R}(e_i,Z)e_i)^{\bot} + (\bar{\nabla}_{mH}
Z)^{\bot} -\nabla_{Z^{\top}}^{\bot}mH 
+\sum_i(\bar{R}(e_i,Z^{\top})e_i)^{\bot}-\tilde{B}(Z^{\bot}).$$
Then the expression of $\mathcal{J}_{\ti{\Omega}}(Z^{\bot})$ 
follows immediately. Now we suppose
 a calibrated extension $M'$ exists with $TM'=EM$ along $M$, and 
$M'$ is minimal. Let $Z'$ and $Z''$  be the  projection of $Z$ onto $TM'$
and ${TM'}^{\bot}$, respectively. Then for 
 $e'_i$ a local o.n.\ frame of $M'$, 
$$div_{M'}(Z')=\sum_i\bar{g}(\bar{\nabla}_{e'_i}Z', e'_i)=
\sum_i-\bar{g}(\bar{\nabla}_{e'_i}Z'', e'_i)
=(m+1)\bar{g}(Z'', H^{M'})=0,$$
where $H^{M'}$ is the mean curvature of $M'$ on $\bar{M}$.
Thus, for $M=\partial O'$,
$$\int_Ma_{Z^{\bot}}dM=\int_{\partial O'}\bar{g}(Z', \nu)=
\int_{O'}div_{M'}(Z')=0. $$
By Proposition 4.1, $Z^{\bot}\in H^1_T(NM)$, and 
supposing $\nabla^{\bot}H=0$, then
Proposition 4.2 yields the equivalence between (14) and the
assumption of $Z^{\bot}$ being an $\Omega$-Jacobi field. As
$\bar{g}(\bar{\nabla}_{\nu}Z,\nu)=0$ holds, then
$c=m\|H\|\, \bar{g}(C_{\ti{\Omega}}(Z^{\bot}),\nu)$.
If $\Omega$ is parallel, then $c=0$ as explained before Lemma 3.1.\qed
\begin{remark}\rm{ If  a normal section
 $W$  is a solution of the $\Omega$-Jacobi operator  
$\mathcal{J}'_{\ti{\Omega,D}}(W)=0$, on a compact domain $D$,
then $\mathcal{J}_{\ti{\Omega}}(W)=c\nu$, where $c=\Psi_{\Omega,D}(W)$
is constant. Supposing
 $W=0$ on a non-empty open set $D'\subset D$,
then  $\mathcal{J}_{\ti{\Omega}}(W)=0$ on $D'$.
 Consequently, $c=0$, and  
$\mathcal{J}_{\ti{\Omega}}(W)=\mathcal{J}'_{\ti{\Omega,D}}(W)=0$ on $D$.
This implies that
$\|\Delta^{\bot} W\|\leq C(\|W\|+ \|\nabla^{\bot}W\|)$, 
 for some constant $C>0$ depending on $D$, $\bar{R}$, $B$, $\Omega$,
$\bar{\nabla}\Omega$, and $\nabla^{\bot}W'_{\alpha}$, where
$W'_{\alpha}$ is a fixed family of o.n. frames of $NM$
defined on a finite cover of $D$ by compact domains. Thus, 
by  Aronszajn's unique continuation theorem for systems of inequalities
of second order (Remark 3 of \cite{Ar}), $W$ must vanish on all $D$. That is, 
$\mathcal{J}'_{\ti{\Omega,D}}$ has uniqueness in the Cauchy problem.
The extra term $C_{\ti{\Omega}}$ can be seen to act in the sense of
distributions.
The coerciveness property
  associated with
$Q_{\Omega}$ still holds on compact domains.
To see this we only have to observe that if $P$ is a bilinear map, then
$|P(W, \nabla^{\bot}_X W')|\leq \|P\|(\lambda\|W\|
+\lambda^{-1}\|\nabla^{\bot}_X W'\|)$ holds for any
 $\lambda>0$, which should be taken  sufficiently
large (see Chap. 8 \cite{GT}).
It follows that a
 Morse index theorem can be stated for submanifolds
with parallel mean curvature and calibrated extended tangent spaces by using
the $\Omega$-Jacobi fields, in a similar way as Simons's version for minimal
submanifolds in \cite{Si} (see also \cite{FriTha}).
}
\end{remark}
If  a calibrated extension $M'$ of $\phi$ exists and 
${\nabla}^{\bot}\nu=0$,
$M$ has parallel mean curvature in $M' $ if and only if  it does in $\bar{M}$. 
If we consider variations
$\phi_t$  with values on $M'$ only, the concepts of
volume-preserving coincide, for $\Omega$ is the volume form on $M'$.
In particular, if  $f\in \mathcal{F}_D$, we  have
\begin{eqnarray*}
{J''}_D(0)(f)&=&\int_D-f\Delta f -(R' +\|B^{\nu} \|^2)f^2 dM
=:I(f,f),\\
&=&\int_D\|\nabla f\|^2-(R' +\|B^{\nu} \|^2)f^2 dM=:q(f), 
\end{eqnarray*}
where 
$R' =\sum_i R'(e_i,\nu,e_i,\nu)=\mathrm{Ricci}'(\nu,\nu)$, with
$\mathrm{Ricci}'$ the Ricci tensor of $M'$, 
and $I$ is the bilinear form defined
in \cite{BdCE}. The immersion into $M'$,
$\phi:D\to M' $, is said to be stable, 
if $I(f,f)\geq 0$ for all $f\in\mathcal{F}_D$.
Considering any section 
$W\in \bar{\mathcal{F}}_{D,\Omega}$
with $W^{\bot}=f\nu$, then $a_W=\bar{g}(W,\nu)=f$,  
and using the Gauss equation for $M'$ as a submanifold of $\bar{M}$, 
we have
\begin{eqnarray}
I_{\Omega}(W,W) 
&=& I(f,f) + m\|H\|\int_Df^2\bar{g}(C_{\ti{\Omega}}(\nu),\nu)dM\nonumber\\
&&+
 \int_Df^2\left( (m+1)\bar{g}(H^{M'},B^{M'}(\nu,\nu))
-\|B^{M'}(\nu,\nu)\|^2\right)dM
\end{eqnarray}
where $B^{M'}$ and $H^{M'}$ stand for the second fundamental form and 
mean curvature of $M'$ in $\bar{M}$, respectively.  We
have used that $B^{M'}(e_i,\nu)=\nabla^{\bot}_{e_i}\nu=0$.
Recall the first eigenvalue
of the twisted Dirichlet problem
 \cite{BaBe} (see also a Euclidean version
\cite{FH}) is given by
$$\lambda_{\mathcal{F}}(D)= \inf\left\{ \frac{q(f)}{\int_D f^2dM}:~~~~
f\in \mathcal{F}_D\right\}.$$
Now we consider the case  $\bar{g}(H^{M'}, B^{M'}(\nu,\nu))=0$ ( for instance,
when $M'$ is minimal, or $B^{M'}(\nu,\nu)=0$). We have
an orthogonal split
$NM=\mathbb{R}\nu \oplus F$ into two parallel subbundles.
If $W\in \mathbb{R}\mathcal{F}_{D,\Omega}$, 
$W=W^{\nu}+W^F= f\nu +W^F$, 
where    $f \in \mathcal{F}_D$ and $W^F\in\mathcal{F}'(F)$, then
\begin{eqnarray}
\lefteqn{I_{\Omega}(W,W) =I(f,f)+ I_{\Omega}(W^F,W^F)
-\int_Df^2(\|B^{M'}(\nu, \nu)\|^2-m\|H\|\bar{g}(C_{\ti{\Omega}}(\nu),\nu))
 dM ~~~~\quad~~~~~~~~}\nonumber\\
&&-2\int_Df\La{(}\sum_i\bar{R}(e_i,\nu,e_i,W^F)
+\sum_{ij}B^{\nu}(e_i,e_j)\bar{g}(B(e_i,e_j),W^F)\La{)}dM.~~~~~\quad\quad\quad
\end{eqnarray}
There are several situations with $C_{\ti{\Omega}}=0$. One is given in 
Proposition 3.2. Another is when $n=2$ and $\bar{\nabla}\Omega=0$. 
In  Lemma 4.4 we will completely
characterize  this condition.
\begin{theorem} We suppose  a  calibrated extension $M'$ of $M$ 
exists satisfying the condition $\bar{g}(H^{M'}, B^{M'}(\nu,\nu))=0$.\\[1mm]
$(1)$ If $\phi:M\to \bar{M}$ is $\Omega$-stable on $D$, 
then
$\phi:M\to M'$ is also stable on $D$ and 
$$\lambda_{\mathcal{F}}(D)\geq \inf_D\left(
\|B^{M'}(\nu,\nu)\|^2-m\|H\|\bar{g}(C_{\ti{\Omega}}(\nu),\nu)\right).$$
$(2)$ If $M'$ is a totally geodesic submanifold of $\bar{M}$
and $C_{\ti{\Omega}}=0$, we have
 for $W\in \mathcal{F}_{D,\mathcal{F}}$ 
$$I_{\Omega}(W,W)=I(f,f)+ I_{\Omega}(W^F,W^F).$$
In the particular case $\bar{g}(\bar{R }(W^F),W^F)\leq 0$, we have
$I_{\Omega}(W,W)\geq  I(f,f)+\int_D\|\nabla^{\bot} W^F\|^2.$
In this case, if $\phi:M\to M'$ is stable, then $\phi:M\to
\bar{M}$ is also $\Omega$-stable.
\end{theorem}
\noindent
\em Proof. \em
It is clear that if $\phi:M\to \bar{M}$ is $\Omega$-stable on $D$ 
then $\phi:M\to M'$
is also stable on $D$. By the assumptions, 
   (15) reads, for $W^{\bot}=f\nu$
with $f\in \mathcal{F}_D$,
\begin{equation}
I_{\Omega}(W,W)= I(f,f)- \int_Df^2\left(\|B^{M'}(\nu,\nu)\|^2-m\|H\|\bar{g}(
C_{\ti{\Omega}}(\nu),\nu) \right)dM.
\end{equation}
Take $f\in \mathcal{F}_D$ an eigenvector of 
$\lambda_{\mathcal{F}}=\lambda_{\mathcal{F}}(D)$ for
the associated twisted Dirichlet problem on $D$, that is,
$-\Delta f- (R'+\|B^{\nu}\|)f=\lambda_{\mathcal{F}}f +\Psi(f)$,
where  $\Psi(f)=|D|^{-1}
\int_D(-\Delta f-(R'+\|B^{\nu}\|)f)$.
Then $I(f,f)=\lambda_{\mathcal{F}}\int_Df^2$. From (17) and the fact that $
I_{\Omega}(W,W)\geq 0$, (1)  follows immediately.
If $M' $ is  totally geodesic, then
$B(e_i,e_j)=B^{\nu}(e_i,e_j)$ takes values on $TM'$, as well as
$\bar{R}(e_i,\nu)e_i=R'(e_i,\nu)e_i$. Thus the last terms  of (16) vanish.
Moreover $\tilde{B}(W^F)=0$. Consequently
$I_{\Omega}(W^F,W^F)=\int_D(\|\nabla^{\bot}W^F\|^2-\bar{g}(\bar{R}(W^F),W^F)) 
dM$. If ~$\bar{g}(\bar{R}(W^F),W^F)\leq 0$,~ then 
$~I_{\Omega}(W^F,W^F)\geq $ 
$\int_D\|\nabla^{\bot}W^F\|^2dM$, which proves (2).\qed \\[2mm]

Barbosa, do Carmo and Eschenburg proved in
\cite{BdCE} that geodesic spheres of space forms are the unique
 stable hypersurfaces
of constant mean curvature. 
The uniqueness is established by showing that
the stability condition implies the hypersurface to be umbilical.
As an immediate consequence  of this result and the  preceding theorem, we have:
\begin{corollary} If a calibrated extension $M'$ exits, and is a space form, 
 if $\phi:M\to\bar{M}$ is $\Omega$-stable  then
$M$ is a geodesic sphere of $M'$.
\end{corollary}
\noindent
In the general case the calibrated extension $M'$ is not a
space form, and so geodesic spheres of $M'$ may  have no constant mean
curvature, nor  be umbilical (the second fundamental form
of geodesic spheres is, up to a sign, the Hessian of the distance 
function to a point),
but umbilical submanifolds may exist. Furthermore,
geodesic $m$-spheres in a Euclidean space or 
in a Euclidean $(m+n)$-sphere with $n\geq 2$  may not be
stable (see Propositions 4.5 and 4.6(2)).
For the case of positive sectional curvature, a more general statement is
the following:
\begin{proposition} Assume $M$ is closed, and $\bar{\nabla}\Omega=0$ or 
$C_{\ti{\Omega}}=0$.  If
$NM$ allows a global unit parallel
section $\nu^F$ orthogonal to $\nu$, and if
$\int_M\sum_i\bar{R}(e_i,\nu^F,e_i,\nu^F)dM> 0$, then $M$
is $\Omega$-unstable.
\end{proposition}
\noindent
\em Proof. \em 
For $\bar{\nabla}\Omega=0$ and $\nabla^{\bot}\nu^F=0$ we also have
$\bar{g}(C_{\Omega}(\nu^F),\nu^F)=0$.
By Proposition 4.1 $\nu^F$  is in the extended domain of
$I_{\Omega}$, and 
$I_{\Omega}(\nu^F,\nu^F)=-\int_M\sum_i\bar{R}(\nu^F,e_i,\nu^F,e_i)+
\sum_{ij}(\bar{g}(B(e_i,e_j),\nu^F))^2<0$. \qed\\[3mm]
We cannot expect Euclidean
spheres to be $\Omega$-stable in $\mathbb{R}^{m+n}$ for any calibration
$\Omega$. We will show in the next proposition how stability depends on 
$\Omega$.
 For any submanifold $M$, consider the tensor
$\xi:\wedge^2NM\to TM^*$ defined by
$$\xi(W,W')(u)=\Omega(W,W',*u)~~~(\mbox{with}~
*\mbox{~the~star~operator~on~}M).$$
\begin{lemma}
The differential operator
$C_{\ti{\Omega}}$ vanishes if and only if $\xi$ vanishes  
and, 
$$
\bar{\nabla}_{W}\Omega(W',d\phi(e_1), \ldots,d\phi(e_m))=
-\bar{\nabla}_{W'}\Omega (W,d\phi(e_1), \ldots,d\phi(e_m))
~~~\forall W,W'\in NM.$$
\end{lemma}
\noindent
\em Proof. \em Fixing a point $p\in M$ and $W_a$ a local o.n. frame of
$NM$ that satisfies $\nabla^{\bot}W_a(p)=0$, we see that
$\bar{g}(C_{\ti{\Omega}}(W_a),W_b)=0$ is equivalent to the last
condition, and taking $W=W_a+fW_b$ where $f$ is any local function, the
condition $\bar{g}(C_{\ti{\Omega}}(W),W)=0$, at $p$, translates into
$ \xi(W_{a},W_{b})(\nabla f)=0$, at $p$.
Since $\nabla f(p)$ is arbitrary, we conclude $\xi=0$.\qed\\[3mm]

\begin{proposition} Suppose $M$ is a
closed pseudo-umbilical submanifold  of $\bar{M}=\mathbb{R}^{m+n}$,
has parallel mean curvature, and  calibrated extended tangent space.
We also suppose $\Omega$ is a calibration on $\mathbb{R}^{m+n}$
that satisfies $\bar{\nabla}_{W}\Omega(W, e_1,\ldots, e_m)=0$ for any 
$W\in NM$. Let $M'$ be 
 the minimal $\Omega$-calibrated  extension  given in Example 2.4, and
 $\lambda_1$ be the first non-zero eigenvalue of $M$
for the closed eigenvalue problem.\\[1mm]
$(1)$ $M$ is stable in $M'$ 
if and only if $\lambda_1\geq m\|H\|^2$. This holds (with equality) when 
$M'=\mathbb{R}^{m+1}$.\\[1mm]
$(2)$ Assume $M$ is totally umbilical in $\bar{M}$, that is, $M'$ is an 
$(m+1)$-Euclidean space
and   $M$ is a Euclidean $m$-sphere of $M'$, and
fix a global parallel basis $W_{\alpha}$ of $TM'^{\bot}=\mathbb{R}^{n-1}$.
Then,   $M$
is $\Omega$-stable in $\bar{M}$  if and only if the  $1$-forms 
$\xi(W_{\alpha},W_{\beta})$ are co-exact, that is, $\xi(W_{\alpha},W_{\beta})=
\delta \omega_{\alpha\beta}$, for some 2-forms $\omega_{\alpha\beta}$ on $M$,
and they satisfy the  inequality
\begin{equation}
 \sum_{\alpha <\beta}
-2m\|H\|\int_M\langle \omega_{\alpha\beta},df_{\alpha}\wedge 
df_{\beta}\rangle dM
\leq \sum_{\alpha}\int_M\|df_{\alpha}\|^2 dM~~~~~
\forall f_{\alpha},f_{\beta}\in C^{\infty}(M)
\end{equation}
where $\langle, \rangle$ denotes the usual Hilbert-Schmidt inner product for
2-forms.
In this case, for each $\alpha, \beta$ 
the following estimates for $\omega_{\alpha\beta}$  holds, viz.
$$\begin{array}{l}
2m\|H\|\left|\int_M \langle \omega_{\alpha\beta},df\wedge dh\rangle
dM \right|\leq \int_M (\|df\|^2+ \|dh\|^2) dM,\\[2mm]
m\|H\|\left|\int_M \omega_{\alpha\beta}(\nabla f, \nabla h)
dM \right|\leq \sqrt{\int_M \|\nabla f\|^2dM}\sqrt{\int_M\|\nabla h\|^2 dM},
\end{array}
$$  
for any smooth functions $f,h$. Furthermore, if $C_{\ti{\Omega}}=0$,
 then $M$ is $\Omega$-stable in $\bar{M}$.
\end{proposition}
\noindent
\em Proof. \em (1) Recall that for $M$ closed, 
$\lambda_1=\inf_{f\in \mathcal{F}_M}(\int_M\|\nabla f\|^2)/(\int_M f^2)$.
Since  $\|B^{\nu}\|^2=m \|H\|^2$ and by Proposition 2.1 in Example 2.4, 
$R'=Ricci'(\nu,\nu)=-\|B^{M'}(\nu,\nu)\|^2=0$, then
$q(f)=\int_M\|\nabla f\|^2dM -m\|H\|\int_Mf^2 dM$, 
and (1) follows from the above  Reighley characterization  of $\lambda_1$.
(2) If $M$ is umbilical in $\bar{M}$, $M$ is a sphere by Proposition 2.1.
From the assumptions, $\bar{g}(C_{\ti{\Omega}}(\nu),\nu)=0$. Then (16) gives us
$I_{\Omega}(W,W)=I(f,f)+ I_{\Omega}(W^F,W^F)$, with
$ I_{\Omega}(W^F,W^F)=\int_M(\|\nabla^{\bot}W^F\|^2 + 
m\|H\|\bar{g}(C_{\ti{\Omega}}(W^F),W^F))dM$. Since $\lambda_1=m\|H\|^2$, 
then  by (1), $I(f,f)\geq 0$, and equality holds for $f$ an
eigenfunction of $\lambda_1$. 
Thus, $M$ is $\Omega$-stable if and only if  $I_{\Omega}(W^F, W^F)\geq 0$.
We take $W_{\alpha}$ a global o.n.\ frame of
parallel sections of $TM'^{\bot}=\mathbb{R}^{n-1}$, and set 
$W^F=\sum_{\alpha}f_{\alpha}W_{\alpha}$, with $f_{\alpha}$ 
arbitrary functions.
Note that $*e_i=(-1)^{i-1}e_1\wedge\ldots\wedge \hat{e}_{i}\wedge 
\ldots\wedge e_m$. 
Then, the $\Omega$-stability condition translates into
\begin{equation}
\int_M \left(\sum_{\alpha}\|\nabla f_{\alpha}\|^2 +\sum_{\alpha\beta}
 m\|H\| f_{\alpha}\xi(W_{\alpha},W_{\beta})(\nabla f_{\beta})\right) dM\geq 0.
\end{equation}
Now we prove the $\xi(W_{\alpha},W_{\beta})$ are co-closed.
If we 
 choose $f_{\alpha}=1$ and $f_{\beta}$ arbitrary, 
and $f_{\gamma}=0$ for $\gamma \neq \alpha, \beta$, we get from (19)
$$\int_M \la{(}\|\nabla f_{\beta}\|^2+m\|H\|\xi(W_{\alpha},W_{\beta})
(\nabla f_{\beta})\la{)}dM
\geq 0.$$
Replacing $f_{\beta}$ by $tf_{\beta}$, with $t>0$ a constant,
 and letting $t\to 0$, we obtain
$\int_M \xi(W_{\alpha},W_{\beta})(\nabla f_{\beta})$ $ \geq 0$, and  again, 
replacing $f_{\beta}$ by $-f_{\beta}$, we obtain equality to zero. 
Since we have  $\xi(W_{\alpha},W_{\beta})(\nabla f_{\beta})=
\langle \xi(W_{\alpha},W_{\beta}),df_{\beta}\rangle$, we conclude
that $\xi(W_{\alpha},W_{\beta})$ are $L^2$-orthogonal
to all exact 1-forms $df_{\beta}$ on $M$. As the Betti numbers of the 
spheres vanish,
by the Hodge decomposition theorem $\xi(W_{\alpha},W_{\beta})$
are co-exact, that is, $\xi(W_{\alpha},W_{\beta})=\delta \omega_{\alpha\beta}$,
for some 2-forms $\omega_{\alpha\beta}$.
Then inequality (19)
 is  equivalent to
\begin{eqnarray*}
\int_M\sum_{\alpha}|df_{\alpha}|^2dM+ \sum_{\alpha\beta}
m\|H\|\int_M\langle \omega_{\alpha\beta},df_{\alpha}\wedge df_{\beta}
\rangle dM \geq 0,
\end{eqnarray*} 
which gives the first inequality of the proposition.
Fixing $\alpha< \beta$ and setting $f=f_{\alpha}$, $h=f_{\beta}$, 
and $f_{\gamma}=0$ for $\gamma\neq \alpha, \beta$,
the above inequality implies
$$\int_M (|df|^2+|dh|^2)dM + 2m\|H\|\int_M \langle \omega_{\alpha\beta},
df\wedge dh\rangle dM \geq 0.$$
If we change $f$ by $-f$, we conclude the second inequality of the 
proposition. 
Note that $\langle \omega_{\alpha\beta},
df\wedge dh\rangle=\omega_{\alpha\beta}(
\nabla f, \nabla h)$.
The last inequality is obtained from the second one
by multiplying $f$ by a constant $t$ and $h$ by $t^{-1}$, with
$t^2=\|\nabla h\|_{L^2}/\|\nabla f\|_{L^2}$.
Finally, if $C_{\ti{\Omega}}=0$, by Lemma 4.4, $\xi(W_{\alpha},
W_{\beta})=0$.\qed \\[3mm]
We  note that for any parallel
calibration  $\Omega$ of $\mathbb{R}^{m+n}$, and any sphere $M$ of a
calibrated vector space $\mathbb{R}^{m+1}$, if we fix $W_{\alpha}$
a constant o.n.\ frame of $TM'^{\bot}=\mathbb{R}^{n-1}$,
then  the  $(m-1)$-forms 
$\hat{\xi}_{\alpha\beta}=\Omega(W_{\alpha},W_{\beta},\cdots)$
are parallel in $\mathbb{R}^{m+n}$.
We may take in previous proposition $\xi(W_\alpha,W_{\beta})=
*\phi^*\hat{\xi}_{\alpha\beta}$, that are obviously co-closed on $M$.
Many well known calibrations in $\mathbb{R}^{m+n}$
satisfy $C_{\Omega}\neq 0$ with $\xi(W_{\alpha},W_{\beta})$
co-closed. 
On the other hand, to investigate if
 inequalities in Proposition 4.5(2) are satisfied or not seems to be
not so easy to determine, as we can see in next remark.
\begin{remark} \rm{ (1)
The associative calibration of $\mathbb{R}^7$ is the 3-form  given by
$$\Omega= \epsilon_*^{123}+\epsilon_*^{145}+ 
\epsilon_*^{167}+\epsilon_*^{246}-\epsilon_*^{257}-
\epsilon_*^{347}-\epsilon_*^{356}$$
where $\epsilon_i$ is the canonical basis of $\mathbb{R}^7$.
We are considering $\mathbb{S}^{2}$ the unit sphere of 
the calibrated subspace spanned by $\epsilon_i$, $i=1,2,3$
 and $W_{\alpha}=\epsilon_{\alpha}$, $\alpha=4,5,6,7$.
Then we have
$ \hat{\xi}_{45}=\hat{\xi}_{67}=\epsilon_*^1=dx^1$, 
$\hat{\xi}_{46}=-\hat{\xi}_{57}=\epsilon_*^2= dx^2$,
$\hat{\xi}_{47}=-\hat{\xi}_{56}=-\epsilon_*^3=-dx^3$.
Consequently, $\omega_{\alpha\beta}=\rho_{\alpha\beta}
\mathrm{Vol}_{\mathbb{S}^2}$ with
${\rho}_{45}={\rho}_{67}=-\phi^1$,
${\rho}_{46}=-{\rho}_{57}=-\phi^2$, 
${\rho}_{47}=-{\rho}_{56}=\phi^3$, 
where $\phi:\mathbb{S}^2\to \mathbb{R}^{3}\subset \mathbb{R}^{7}$ 
is the inclusion map.
Let us suppose that $\phi$ is $\Omega$-stable. Then $\omega_{45}$
should satisfy the last inequality of Proposition 4.5, i.e.\ 
for any functions $f,h:\mathbb{S}^2\to \mathbb{R}$
$$2 \left|\int_{\mathbb{S}^2} \phi^1  
\mathrm{Vol}_{\mathbb{S}^2}(\nabla f,\nabla h)dM \right|\leq 
\sqrt{\int_{\mathbb{S}^2}\|\nabla f\|^2dM}
\sqrt{\int_{\mathbb{S}^2}\|\nabla h\|^2dM}$$
We now use the stereographic projection
$\sigma:\mathbb{R}^2\to \mathbb{S}^2\subset \mathbb{R}^3$, 
$\sigma(w)=\left(\frac{(|w|^2-1)}{(|w|^2+1)}, \frac{2w}{(|w|^2+1)}\right),$
 that is a conformal map.
We denote by $\mathrm{Vol}_0$ the Euclidean volume element of $\mathbb{R}^2$,
$J$ the canonical complex structure,
and by $\nabla^0f$ the Euclidean gradient  
for a function defined on $\mathbb{R}^2$. 
Then the above
inequality is equivalent to
$$2\left|\int_{\mathbb{R}^2}\left(\frac{|w|^2-1}{|w|^2+1}\right)
\mathrm{Vol}_{0}(\nabla^0f,\nabla^0h)dw \right| \leq 
\sqrt{\int_{\mathbb{R}^2}|\nabla^0f|^2dw}
\sqrt{\int_{\mathbb{R}^2}|\nabla^0h|^2dw}
$$
for functions $f, h :\mathbb{R}^2\to \mathbb{R}$, that we take
with compact support in an annulus $D:=\{w: 0\leq  R_1\leq |w|\leq R_2\}$.
We choose $R_1$ sufficiently large so that ${(|w|^2-1)}/{(|w|^2+1)}\geq
\frac{1}{2}+ \delta$, where $0<\delta<1/2$ is a constant.
Since $g_0(J\nabla^0 f,\nabla^0 h)=\mathrm{Vol}_0(\nabla^0 f,\nabla^0 h)$, and
if this is $\geq 0$, 
from preceding inequality
we have
$2(\frac{1}{2}+\delta)\left|\int_{D}
g_0(J\nabla^0f, \nabla^0h) dw\right| \leq
|\nabla^0f|_{L^2}|\nabla^0h|_{L^2}$. 
A pair of functions $(h,h')$ on $D$ defines 
a holomorphic map in $\mathbb{C}$ if and only if
 $-J\nabla^0 h'(w)=\nabla^0 h(w)$. Thus,
we maximize $g_0(J\nabla^0f,\nabla h)$ by taking $f=-h'$
for such pair of conjugate harmonic maps, giving $
\mathrm{Vol}_0(\nabla^0 f,\nabla^0 h)=
\|\nabla^0 f\|
\|\nabla^0 h\|=\|\nabla^0 f\|^2$. 
This would give a contradiction in the previous inequality. As a matter 
of fact, we cannot choose such a pair of functions, 
because nonconstant harmonic maps cannot vanish in all $\partial D$.
We also recall that a
harmonic function $f$ on $\mathbb{R}^2$ with $L^2$ derivative
defines a $L^2$-harmonic one-form $df$, and so, by a result of
Yau ( \cite{Yau}, Theorem 6), $f$ must be constant. 
The question is to know how far is 
an holomorphic  map $(\varphi,\varphi')$ 
on $D$ from a pair of  functions $(f,h)$  vanishing in $\partial D$. 
This can be measured by  $|\mathrm{Vol}_0(\nabla^0f,\nabla h)|$ 
as we have described. A significant distance between these two set
of pairs of functions could indicate that (18) holds.
We also observe that $\phi^k$ are  $\lambda_1$-eigenfunctions 
with $\lambda_1=2$  the first nonzero eigenvalue of $\mathbb{S}^2$,
and they satisfy $\int_{\mathbb{S}^2}(\phi^k)^2dM=\frac{1}{3}|\mathbb{S}^2|$.
Furthermore, the inequality   (18) holds 
if we take $f_{\alpha}$ and $f_{\beta}$ any  $\lambda_i$-eigenfunctions,
with $i=1$ or $2$, giving either equality, or zero in the l.h.s.
The  fact that $\omega_{\alpha\beta}$ is defined using the
$\lambda_1$-eigenfunctions suggests us that a
 proof of the $\Omega$-stability of the
2-sphere should be  related to some  inequalities
derived from spectral theory, and this will be the
subject of future work.
\\[2mm]
(2) Let us now consider the K\"{a}hler calibration of $\mathbb{R}^6$ 
given by the 4-form
$$\Omega=\frac{1}{2}(\epsilon_*^{12}+\epsilon_*^{34}+\epsilon_*^{56})
\wedge (\epsilon_*^{12}+\epsilon_*^{34}+\epsilon_*^{56})
=\epsilon_*^{1234}+\epsilon_*^{1256}+\epsilon_*^{3456}$$
where $\epsilon_i$ is the canonical basis of $\mathbb{R}^7$,
and $\phi=(\phi^1,\ldots,\phi^4)
:\mathbb{S}^3\to \mathbb{R}^4$ denotes the inclusion map of the
$3$-sphere of
 the calibrated subspace $\mathbb{R}^4$
spanned by $\epsilon_i$, $i=1,2,3, 4$.
Then we are taking
$\hat{\xi}_{56}=\epsilon_*^{12}+\epsilon_*^{34}=
dx^1\wedge dx^2+dx^3\wedge dx^4$. 
The $\Omega$-stability condition is equivalent to the inequality
$$3\left|\int_{\mathbb{S}^3}\langle \omega_{56},df_5\wedge df_6
\rangle dM\right|
\leq\sqrt{\int_{\mathbb{S}^3}\|\nabla f_5\|^2dM}
\sqrt{\int_{\mathbb{S}^3}\|\nabla f_6\|^2dM}$$
to be valid for all smooth maps $f_5,f_6$ on $\mathbb{S}^3$.
We have
$*\omega_{56}=\phi^1d\phi^2+\phi^3d\phi^4$,
and
$\langle \omega_{56}, df_5\wedge df_6\rangle dM=
df_5\wedge df_6 \wedge *\omega_{56}.$
If we take $f_{\alpha}$ one of the components $\phi^i$ ( that
are $\lambda_1$-eigenfunctions, with $\lambda_1=3$ the first non-zero
eigenvalue of $\mathbb{S}^3$), we can verify, using spherical
coordinates, that the previous inequality holds, with equality
in some cases. Once more, a proof for  stability seems
to be related to new spectral inequalities as in preceding case (1).
}
\end{remark}

Next we  obtain a  uniqueness theorem which extends the case $n=1$ \cite{BdC}:
\begin{theorem} Assume that $\bar{M}=\mathbb{R}^{m+n}$, 
$\bar{\nabla}\Omega=0$ (or $C_{\ti{\Omega}}=0$), and 
$M$ is a closed submanifold with parallel mean curvature and
calibrated extended tangent space.
Consider the height functions  
$$h=\bar{g}(\phi,\nu)~~~~\mbox{and}~~~~
S=\sum_{ij}\bar{g}(\phi, (B(e_i,e_j))^F)B^{\nu}(e_i,e_j).$$
If $\phi:M\to \bar{M}$ is $\Omega$-stable and
$\int_M S(2+h\|H\|)dM\leq 0$,
then $\phi$ is pseudo-umbilical and a minimal calibrated extension
$M'$ of $M$ exists with $R'=-\|B^{M'}(\nu,\nu)\|^2=0$ and $S=0$. 
Furthermore, if $NM$ is a 
trivial bundle, then  $M'=\mathbb{R}^{m+1}$, 
 $M$ is a Euclidean sphere and $\Omega$
satisfies Proposition 4.5(2).
\end{theorem}
\noindent
\em Proof. \em
If $\phi$ is pseudo-umbilical, then by Proposition 2.1 a (minimal) calibrated
extension $M'$ exists, satisfying
$R'= -\|B^{M'}(\nu,\nu)\|^2=0$, and $S=-\bar{g}(B^{M'}(\nu,\nu),\phi)H$ $=0$.
We follow \cite{BdC}.  Let $\bar{X}_x=x$ be the 
position vector field
in $\mathbb{R}^{m+n}$. Using the well-known expression
$m\bar{g}(H,\phi)= \mathrm{div}(\phi^{\top})-\frac{1}{2}tr_gL_{\bar{X}}\bar{g}$,
where $\phi^{\top}(x)$ is the projection of $\phi(x)$ onto $T_xM$,
and integrating over $M$, we have 
\begin{equation}\label{estaF}
\int_M (\|H\| h +1) dM = 0.
\end{equation}
That is, $f=(\|H\| h +1)\in \mathcal{F}_M$. Now,  
$dh(e_i)=\bar{g}(\phi,d\nu(e_i))= \sum_j-\bar{g}(\phi, e_j)B^{\nu}(e_j,e_i)$.
By applying  Coddazzi's equation we have $\sum_i\nabla_{e_i}B^{\nu}(e_i,e_j)=
m\bar{g}(\nabla^{\bot}_{e_j}H,\nu )=0$. 
We may assume   at a given point $\nabla_{e_i}e_j=0$. 
Then at that point
\begin{eqnarray*}
\Delta h &=& \sum_{ij} -\delta_{ij}B^{\nu}(e_i,e_j)
-\bar{g}(\phi, B(e_i,e_j))B^{\nu}(e_i,e_j)
-\bar{g}(\phi, e_j)\nabla_{e_i}B^{\nu}(e_j,e_i)\\
&=& -m\|H\| -\bar{g}(\phi, B(e_i,e_j))B^{\nu}(e_i,e_j)=-m\|H\|
-h\|B^{\nu}\|^2- S,
\end{eqnarray*}
and integration over $M$ give us the equality
\begin{equation}\label{outraigualdade}
\int_M h \|B^{\nu}\|^2dM=-\int_M(S+m\|H\|)dM.
\end{equation}
The stability condition  applied to $W=f\nu$, and the fact that 
$\bar{g}(C_{\ti{\Omega}}(\nu),\nu)=0$, by
assumptions on $\Omega$ and $\nu$, implies
$\int_M -f\Delta f - f^{2}\|B^{\nu}\|^2\geq 0$, that is,
\begin{eqnarray*}
&&\int_M\La{(}\|H\|^2h^2\|B^{\nu}\|^2 + m\|H\|^3h +\|H\|^2 h S 
+\|H\|h\|B^{\nu}\|^2 +m\|H\|^2 + \|H\| S\La{)} dM ~~~~~~~ \\
 &&\geq\int_M(\|B^{\nu}\|^2|H\|^2 h^2+ 2 \|B^{\nu}\|^2\|H\|h 
+ \|B^{\nu}\|^2) dM, 
\quad\quad\quad
\end{eqnarray*}
and using the above equalities (\ref{estaF})(\ref{outraigualdade}),  
we get the simplified inequality
$$ \int_M \|H\| S(\|H\|h +1)dM \geq \int_M
\|B^{\nu}\|^2(\|H\| h+1) dM=\int_M \|B^{\nu}\|^2 -\|H\|(S+m\|H\|)dM,$$
that is,  $\int_M \|H\|S(\|H\| h+2)dM \geq \int_M (\|B^{\nu}\|^2-m\|H\|^2)dM$.
By assumption, and the fact that $\|B^{\nu}\|^2\geq m\|H\|^2$,  we get
$\|B^{\nu}\|^2=m\|H\|^2$, which proves $M$ is pseudo-umbilical. Thus,
$S=0$.
If $NM$ is spanned by a global orthonormal system of $n$ parallel  sections 
$\{\nu, W_{\alpha}\}$,
then $\bar{g}(C_{\ti{\Omega}}(W_{\alpha}),W_{\alpha})=0$ and the
  $\Omega$-stability implies 
$$I_{\Omega}(W_{\alpha},W_{\alpha})=-\int_M\sum_{ij}\bar{g}(B(e_i,e_j),
W_{\alpha})^2\geq 0,$$ 
that is, $B^{M'}(e_i,e_j)=0$. Thus $\phi$ is
totally umbilical and, by Proposition 2.1, $M$ is a sphere on
$M'=\mathbb{R}^{m+1}$
and $C_{\ti{\Omega}}$ satisfies Proposition 4.5(2).
\qed\\[3mm]

Now we specialize on the case that the calibration is $\Omega_{\pi}$, 
defined by a Riemannian fibration $\pi:\bar{M}\to N$ of totally 
geodesic fibres.  We recall that the Riemannian submersions of
the unit Euclidean spheres $\mathbb{S}^{m+n}$, with totally
geodesic and connected fibres, were classified by Escobales and 
Ranjan \cite{Es,Ra},
 and define the Hopf fibrations
of the spheres. Among these fibrations, the ones that have
fibres of dimension $\geq 3$  are  
the Hopf fibrations of $\mathbb{S}^7$ with fibre
$\mathbb{S}^3$, of $\mathbb{S}^{4k+3}$, for $k\geq 2$,  with fibre 
$\mathbb{S}^3$,
and  of $\mathbb{S}^{15}$ with fibre $\mathbb{S}^7$, i.e.
\begin{equation}
 \mathbb{S}^3\hookrightarrow\mathbb{S}^7
\to \mathbb{S}^4(\sm{\frac{1}{2}}),  ~~~~~\mathbb{S}^3\hookrightarrow
\mathbb{S}^{4k+3}\to\mathbb{HP}^k,~~~~~
\mathbb{S}^7\hookrightarrow \mathbb{S}^{15}\to
\mathbb{S}^8(\sm{\frac{1}{2}}),
\end{equation}
respectively, 
where $\mathbb{HP}^k$ is the quaternionic projective space of sectional
curvature $K$ with $1\leq K\leq 4$, and $\mathbb{S}^4(\sm{\frac{1}{2}})$
and $\mathbb{S}^8(\sm{\frac{1}{2}})$ are spheres of curvature $4$. 
Fibrations of $\mathbb{H}^{m+2}$
by totally geodesic  hypersurfaces (and so by $(m+1)$-dimensional
hyperbolic spaces) were described by Ferus \cite{Fe}, and they
arise as the nullity foliation of a suitable isometric immersion
of $\mathbb{H}^{m+2}$ into $\mathbb{H}^{m+3}$ without umbilics.

\begin{proposition} $(1)$ Any Euclidean $m$-dimensional sphere of 
 an $(m+1)$-dimensional vector  subspace $ E$
is $\Omega_{\pi}$-stable in $\mathbb{R}^{m+n}$ for any fibration 
in $\mathbb{R}^{m+n}$  with $E$ as a fibre.
\\[2mm]
$(2)$ Let $M^m$ be a geodesic sphere of
$ \mathbb{S}^{m+1}$. Immersing $ \mathbb{S}^{m+1}$ as a totally geodesic
fibre in $\mathbb{S}^{m+n}$, where $m, n$ are such that
 $\pi:\mathbb{S}^{m+n} \to N$ is
 one of the Hopf fibrations  given in (22), defines an
 $\Omega_{\pi}$-unstable  immersion $\phi:M^m\to \mathbb{S}^{m+n}$
 with parallel mean curvature.
\\[2mm]
$(3)$ If  $M^m$ is a geodesic sphere of a 
hyperbolic space $\mathbb{H}^{m+1}$, with $\mathbb{H}^{m+1}$  immersed
as a fibre of a Riemannian fibration $\pi:\mathbb{H}^{m+n}\to N$ of
the $(m+n)$-dimensional hyperbolic space,  with  $n\geq 2$, and 
 by totally geodesic fibres, then the corresponding 
immersion $\phi:M^m\to \mathbb{H}^{m+n}$ is an  $\Omega_{\pi}$-stable
immersion with parallel mean curvature.
\end{proposition}
\noindent
\em Proof. \em
(1) and (3) are immediate consequences of Theorem 4.1(2)
and Proposition 3.2. 
Now we prove (2). In \cite{BdCE} it is proved that
$I(f,f)\geq 0$ for all $f\in \mathcal{F}_M$. 
On the other hand, the normal bundle $N(\mathbb{S}^{m+1})$
of $\mathbb{S}^{m+1}$ in $\mathbb{S}^{m+n}$ is a trivial bundle, spanned
by $n-1$ unit parallel vector fields $V_1, \ldots, V_{n-1}$. Their restrictions
to $M$ are  parallel along $M$, and they span  the parallel subbundle 
$F$ of $NM$. By  Proposition 4.4 with $C_{\ti{\Omega}}=0$ (see Proposition
3.2), $M$ is $\Omega$-unstable.
\qed\\[3mm]
\noindent
{\bf Acknowledgements.} The author is indebted to
Pedro Freitas for his help on handling  the
twisted eigenvalue problem and related spaces of functions.


\begin{thebibliography}{99}
{\small 
\bibitem{Ar} Aronszajn, N.: {\sl A unique continuation theorem for solutions 
of elliptic partial differential equations or inequalities of second order}.  
J.\ Math.\ Pures Appl.\ (9) {\bf 36}  (1957), 235--249. \\[-7mm]
\bibitem{BdC} Barbosa, J.L., do Carmo, M.:
{\sl Stability of minimal surfaces and eigenvalues of the Laplacian
}. Math.\ Z.\  {\bf 173}  (1980), no.\ 1, 13--28.\\[-7mm]
\bibitem{BdCE} Barbosa, J.L., do Carmo, M., Eschenburg, J.:
{\sl Stability of hypersurfaces of constant mean curvature in 
Riemannian manifolds}.  Math.\ Z. {\bf 197}  (1988),  no.\ 1, 123--138.\\[-7mm]
\bibitem{BaBe} Barbosa, J.L., B\'{e}rard, P.: {\sl  Eigenvalue and ``twisted'' 
eigenvalue problems, applications to CMC surfaces}.  J.\ Math.\ Pures 
Appl. (9) {\bf 79} (2000),  no.\ 5, 427--450.\\[-7mm]
\bibitem{Chen} Chen, B.-Y.:{\sl Geometry of Submanifolds.} Pure and
 Applied Mathematics, No. 22. Marcel Dekker, Inc., New York, 1973.\\[-7mm]
\bibitem{CY} Chen, B.-Y., Yano, K.: {\sl Integral Formulas for submanifolds and
their applications.} J.\ Diff.\ Geom.\ {\bf 5} (1971), 467-477.\\[-7mm]
\bibitem{ChNo} Cheng, Q-M, Nonaka, K.: {\sl Complete submanifolds in Euclidean 
spaces with parallel mean curvature vector.} Manuscripta Math. {\bf 105}
(2001), 353-366.\\[-7mm]
\bibitem{DF1} Duzaar, F., Fuchs, M.:{\sl On the existence of integral
currents with prescribed mean curvature vector.} Manuscripta Math. {\bf 67}
(1990), 41-67.\\[-7mm]
\bibitem{DF2} Duzaar, F., Fuchs, M.: {\sl On  integral
currents with constant mean curvature vector.} Rend.\ Sem.\ Mat.\ Univ.\
Padova, {\bf 85}
(1991), 79-103.\\[-7mm]
\bibitem{Es} Escobales, R.H.: {\sl Riemannian submersions with totally 
geodesic fibres}. J.\ Differential Geom.\ {\bf 10} (1975), 253-276.\\[-7mm]
\bibitem{Fe} Ferus, D.: {\sl On isometric immersions between hyperbolic
spaces}. Math.\ Ann.\ {\bf 205} (1973), 193--200.\\[-7mm]
\bibitem{FH} Freitas, P., Henrot, A.: {\sl On the first twisted Dirichlet 
eigenvalue}.  Comm.\ Anal.\ Geom.\ {\bf 12} (2004),  no.\ 5, 1083--1103. 
\\[-7mm]
\bibitem{FriTha} Frid, H.; Thayer, F.J., {\sl An abstract version 
of the Morse index theorem and its application to hypersurfaces of 
constant mean curvature}.  Bol.\ Soc.\ Brasil Mat.\ (N.S.) {\bf 20} (1990),  
no. 2, 59--68.\\[-7mm] 
\bibitem{GT} Gilbarg, D.\,  Trudinger, N.S.:{\sl  Elliptic partial 
differential equations of second order}, Reprint of the 1998 edition. 
Classics in Mathematics. Springer-Verlag, Berlin, 2001.\\[-7mm]
\bibitem{Gu1} Gulliver, R.:{\sl  Existence of surfaces with prescribed mean
curvature vector}, Math.\ Z.\ {\bf 131} (1973), 117--140.\\[-7mm]
\bibitem{Gu2} Gulliver, R.: {\sl Necessary conditions for submanifolds and 
currents with prescribed mean curvature vector}, 
in Seminar on Minimal Submanifolds, Enrico Bombieri, ed.,
Ann.\ of Math.\ Studies 103, Princeton Univ. Press, 1983, pp 225--242.\\[-7mm]
\bibitem{HL} Harvey, R., Lawson, H.B. Jr.:{\sl  Calibrated geometries}.  
Acta Math.\ {\bf 148}  (1982), 47--157.\\[-7mm]
\bibitem{LS} Li, G., Salavessa, I.M.C.: {\sl Bernstein-Heinz-Chern results
in calibrated manifolds}.  Rev.\ Mat.\ Iberoamericana {\bf 26}(2) (2010), 651--692.
\\[-7mm]
\bibitem{Mo} Morgan, F.:{\sl Perimeter-minimizing curves and surfaces
in $\mathbb{R}^n$ enclosing prescribed multi-volume.} Asian J.\ Math.\
{\bf 4} no.\ 2 (2000), 373--383. \\[-7mm]
\bibitem{Ra} Ranjan, A.:{\sl Riemannian submersions of spheres
with totally geodesic fibres}. Osaka J.\ Math. {\bf 22} (1985), 243--260.
\\[-7mm]
\bibitem{RiRo} Ritor\'{e}, M., Ros, A.: {\sl Stable constant mean curvature 
tori and the isoperimetric problem in three space forms}.  
Comment.\ Math.\ Helv.\ {\bf 67} (1992),  no.\ 2, 293--305.\\[-7mm]
\bibitem{Si} Simons, J.: {\sl Minimal varieties in Riemannian manifolds}.
Ann.\ Math. {\bf 88} (1968) 62--105.\\[-7mm]
\bibitem{Smale} Smale, S.: {\sl On the Morse index theorem}.  
J.\ Math.\ Mech.\  {\bf 14} (1965), 1049--1055.\\[-7mm]
\bibitem{Yau} Yau, S.T.: {\sl Some function-theoretic properties of complete 
Riemannian manifolds and their applications to geometry}. 
Indiana Math.\ J.\ {\bf 25} (1976), 659-670.
}
\end{thebibliography}
\end{document}